\documentclass[11pt]{amsart}
\usepackage{amsmath}
\usepackage{amssymb}
\usepackage{amsfonts}
\usepackage{amsthm}
\usepackage{enumerate}
\usepackage[mathscr]{eucal}
\usepackage{graphicx}
\usepackage[bookmarksnumbered, bookmarksopen, colorlinks, citecolor=blue, linkcolor=blue]{hyperref}

\usepackage{eufrak}
\usepackage{color}
\usepackage{tikz-cd}
\usepackage{enumerate}
\usepackage{cleveref}

\usetikzlibrary{decorations.pathreplacing}

\newtheorem{theorem}{Theorem}[section]
\newtheorem{lemma}[theorem]{Lemma}

\theoremstyle{definition}

\newtheorem{proposition}[theorem]{Proposition}

\theoremstyle{remark}
\newtheorem{remark}[theorem]{Remark}

\newtheorem*{restatedtheorem}{\restatedtheoremnumber}

\newtheorem*{restatedlemma}{\restatedlemmanumber}

\newtheorem*{restatedproposition}{\restatedpropositionnumber}
\newenvironment{restateproposition}[1]{%
  \begin{samepage}\def\restatedpropositionnumber{Proposition~\ref*{#1}}%
  \begin{restatedproposition}%
}{\end{restatedproposition}\end{samepage}\nopagebreak[4]}

\numberwithin{equation}{section}

\def\intslash{\rlap{\kern  .32em $\mspace {.5mu}\backslash$ }\int}
\def\qsl{{\rlap{\kern  .32em $\mspace {.5mu}\backslash$ }\int_{Q_x}}}

\newcommand{\Be}{\begin{equation}}
\newcommand{\Ee}{\end{equation}}
\newcommand{\Bes}{\begin{equation*}}
\newcommand{\Ees}{\end{equation*}}
\newcommand{\Bsp}{\begin{split}}
\newcommand{\Esp}{\end{split}}
\newcommand{\Bm}{\begin{multline}}
\newcommand{\Em}{\end{multline}}
\newcommand{\Bea}{\begin{eqnarray}}
\newcommand{\Eea}{\end{eqnarray}}
\newcommand{\Beas}{\begin{eqnarray*}}
\newcommand{\Eeas}{\end{eqnarray*}}
\newcommand{\Benu}{\begin{enumerate}}
\newcommand{\Eenu}{\end{enumerate}}
\newcommand{\Bi}{\begin{itemize}}
\newcommand{\Ei}{\end{itemize}}

\begin{document}

\title[An improved estimate for cubic \(L^8\) restriction]{An improved estimate for discrete \(L^8\) restriction on the cubic curve}

\author{Xudong Lai}
\address{Xudong Lai: Institute for Advanced Study in Mathematics\\
Harbin Institute of Technology\\
Harbin
150001\\
China;
Zhengzhou Advanced Research Institute\\
Harbin Institute of Technology\\
Zhengzhou
450000\\
China}
\email{xudonglai@hit.edu.cn\ xudonglai@mail.bnu.edu.cn}

\author{Hongxu Wu}
\address{Hongxu Wu: Institute for Advanced Study in Mathematics, Harbin Institute of Technology, Harbin, 150001, China}
\email{2022110278@stu.hit.edu.cn}



\begin{abstract}
We prove an improved estimate for discrete \(L^8\) restriction on the cubic curve \((n,n^3)\). More precisely, for every \(\varepsilon>0\), the following eighth moment estimate holds
\[
\int_{[0,1]^2}
 \left|\sum_{n=1}^N b_n\exp\!\bigl(2\pi i(nx+n^3y)\bigr)\right|^8dxdy\ll_\varepsilon N^{1/2-1/47500+\varepsilon}
 \left(\sum_{n=1}^N|b_n|^2\right)^4.
\]
\end{abstract}
\maketitle
\section{Introduction and main result}\label{sec:introduction}
Discrete Fourier restriction on the cubic curve \((n,n^3)\) arises in the study of periodic Strichartz estimates for the Korteweg--de Vries equation (see \cite{bourgainkdv}).
Let \(N\) be a positive integer, let \(b_1,\ldots,b_N\) be complex numbers, and let \(r>0\). Write
\[
 I_r(b)=\int_{[0,1]^2}
 \left|\sum_{n=1}^N b_n\exp\!\bigl(2\pi i(nx+n^3y)\bigr)\right|^r\,dx\,dy .
\]

Bourgain \cite[Eq.~(8.37)]{bourgainkdv} established the \(L^6\) restriction estimate with an arbitrarily small power of \(N\). For \(r\ge8\), the conjectured bound with an \(\varepsilon\)-loss, in moment form, is
\[
 I_r(b)\ll_{r,\varepsilon} N^{r/2-4+\varepsilon}
       \left(\sum_{n=1}^N|b_n|^2\right)^{r/2} .
\]
The case \(r=8\) is the  critical point where the conjectured power of \(N\), apart from an \(\varepsilon\) loss, is zero.

Hu and Li \cite[Theorem 1.1]{huli} proved the conjectured bound for \(r\ge14\). The first author and Ding \cite[Theorem 1.1 and Corollary 1.4]{laiding} used estimates for the full moment curve from Bourgain, Demeter, and Guth \cite{bourgain-demeter-guth-vmvt} to establish a discrete restriction theorem for curves with monomial coordinates. For \((n,n^3)\), their theorem extended the range to \(r\ge12\), giving in particular
\[
 I_{12}(b)\ll_\varepsilon N^{2+\varepsilon}
       \left(\sum_{n=1}^N|b_n|^2\right)^6 .
\]
Hughes and Wooley \cite[Theorem 1.2]{hugheswooley} reached \(r=10\), proving
\begin{equation}\label{eq:weighted-ten}
 I_{10}(b)\ll_\varepsilon N^{1+\varepsilon}
       \left(\sum_{n=1}^N|b_n|^2\right)^5 .
\end{equation}
There are also several work  treating related restriction problems beyond the cubic curve, we refer to see \cite{huli-schrodinger, huli-fifth, li-mean-value, li-stein-decoupling, li-stein-2025}.

We also recover the sixth-moment bound on positive support by combining the triple representation bound \eqref{eq:triples}, proved in Section~\ref{sec:support}, with Cauchy--Schwarz. This gives
\begin{equation}\label{eq:weighted-six}
 \begin{aligned}
 I_6(b)
 &=\sum_{(s,v)\in\mathbb Z^2}
 \Big|\sum_{\substack{1\le a_1,a_2,a_3\le N\\
 a_1+a_2+a_3=s\\a_1^3+a_2^3+a_3^3=v}}
 b_{a_1}b_{a_2}b_{a_3}\Big|^2\\
 &\ll_\varepsilon N^\varepsilon
 \sum_{a_1,a_2,a_3=1}^N|b_{a_1}|^2|b_{a_2}|^2|b_{a_3}|^2
 =N^\varepsilon\left(\sum_{n=1}^N|b_n|^2\right)^3 .
 \end{aligned}
\end{equation}
Interpolating between \eqref{eq:weighted-six} and \eqref{eq:weighted-ten} by H\"older's inequality, we obtain
\begin{equation}\label{eq:interpolated-eight}
 I_8(b)\le I_6(b)^{1/2}I_{10}(b)^{1/2}
 \ll_\varepsilon N^{1/2+\varepsilon}
       \left(\sum_{n=1}^N|b_n|^2\right)^4 .
\end{equation}

In this paper, we improve this eighth-moment bound by a fixed power of \(N\).

\begin{theorem}\label{th:main}
For every \(\varepsilon>0\), there is a constant \(C_\varepsilon>0\) such that, for every positive integer \(N\) and all complex numbers \(b_1,\ldots,b_N\),
\[
 I_8(b)\le C_\varepsilon N^{1/2-1/47500+\varepsilon}
       \left(\sum_{n=1}^N|b_n|^2\right)^4 .
\]
\end{theorem}

Taking eighth roots in Theorem~\ref{th:main} gives the following \(L^8\) norm bound:
\[
 \left\|\sum_{n=1}^N b_n\exp\!\bigl(2\pi i(nx+n^3y)\bigr)
 \right\|_{L^8([0,1]^2)}
 \ll_\varepsilon N^{1/16-1/380000+\varepsilon}
       \left(\sum_{n=1}^N|b_n|^2\right)^{1/2}.
\]
The main strategy in this paper is as follows. We first work with indicator functions. If the eighth moment is large, we will obtain upper and lower bounds for the size of its support. We then find a large set of integer chord points in a short horizontal interval for which the difference set is small. To bound the difference set from below, we combine finite differences of integer slopes with divisibility estimates and counts of additive cubes. This argument applies to general monomial curves. For the cubic curve, taking five differences contradicts the bounds obtained from a large eighth moment. Finally we pass from indicator functions to arbitrary complex coefficients.

\subsection*{Outline of the paper} In Section~\ref{sec:core} we prove our main theorem from five propositions. We first reduce to indicator functions, then control the support size, and finally compare upper and lower bounds for the difference set of a local chord set. The five propositions are proved in Sections~\ref{sec:support}--\ref{sec:growth}.

\subsection*{Notation} Throughout this paper, \(\mathbf1_A\) denotes the indicator function of \(A\). All solution counts refer to ordered tuples. For finite sets \(U,V\subseteq\mathbb Z^2\), an additive quadruple in \(U\) is a tuple \((P_1,P_2,P_3,P_4)\in U^4\) satisfying \(P_1+P_2=P_3+P_4\). We write
\(U\pm V=\{u\pm v:u\in U,\ v\in V\}\), with matching signs, and \(-U=\{-u:u\in U\}\). The set \(jU\) consists of sums of \(j\) elements of \(U\), with repetitions allowed. For \(z\in\mathbb Z^2\), \(pz\) denotes coordinatewise multiplication by \(p\). For a prime \(p\), \(\mathbb Z/p\mathbb Z\) is the field of residue classes modulo \(p\), and congruences of points are interpreted coordinatewise.
For a positive integer \(m\), let \(\tau(m)\) be the number of its positive divisors. We use the standard divisor bound
\(\tau(m)\ll_\eta m^\eta\) for every \(\eta>0\) \cite[Chapter XVIII]{hardy-wright}. The logarithm is natural.
\(\theta\) is fixed and \(N\) tends to infinity. In an upper bound, a factor \(N^{o(1)}\) is bounded by \(C_{\theta,\eta}N^\eta\) for every \(\eta>0\). In a lower bound, \(N^{-o(1)}\) means the reciprocal of such a factor. The constants are uniform in the sets, horizontal intervals, and chosen primes.

\section{Core Propositions}\label{sec:core}
We collect the estimates leading to Theorem~\ref{th:main}. The first proposition reduces the problem to indicator functions. The remaining propositions describe the consequences of a large eighth moment and give the difference-set estimate that rules them out.

For \(A\subseteq\{1,\ldots,N\}\), let \(\mathbf1_A\) denote its indicator function.

\begin{proposition}\label{prop:weighted}
Fix \(N\ge 1\). If \(D>0\) satisfies \(I_8(\mathbf1_A)\le D|A|^4\) for every \(A\subseteq\{1,\ldots,N\}\), then for all complex coefficients \(b_1,\ldots,b_N\),
\begin{equation}\label{eq:weighted-reduction}
 I_8(b)\le D(1+\log N)^4
       \left(\sum_{n=1}^N|b_n|^2\right)^4 .
\end{equation}
\end{proposition}

We prove Proposition~\ref{prop:weighted} in Section~\ref{sec:growth}. It reduces the main proof to the indicator estimate.
It therefore suffices to show that, for each fixed \(1/2-1/47500<\theta<1/2\),
\begin{equation}\label{eq:indicator-target}
 I_8(\mathbf1_A)\le C_\theta N^\theta|A|^4
 \qquad(N\ge1,\ A\subseteq\{1,\ldots,N\}).
\end{equation}

For indicator functions, orthogonality turns the eighth moment into a count of ordered solutions:
\begin{equation}\label{eq:eight-solutions}
 I_8(\mathbf1_A)=
 \#\left\{(a_1,\ldots,a_4,c_1,\ldots,c_4)\in A^8:
 \begin{array}{l}
 \sum_{j=1}^4a_j=\sum_{j=1}^4c_j,\\
 \sum_{j=1}^4a_j^3=\sum_{j=1}^4c_j^3
 \end{array}\right\}.
\end{equation}
We first ask what the condition \(I_8(\mathbf1_A)\ge N^\theta|A|^4\) tells us about \(|A|\). The next proposition gives both an upper bound and a lower bound.

\begin{proposition}\label{th:support}
Let \(13/27<\theta<1/2\), and let \(A\) be a nonempty subset of \(\{1,\ldots,N\}\).
If \(I_8(\mathbf1_A)\ge N^\theta|A|^4\), then
\[
 N^{\theta-o(1)}\le |A|\le N^{1-\theta+o(1)} .
\]
\end{proposition}

\begin{remark}
Proposition \ref{th:support} gives
\[
 N^{1/2-(1/2-\theta)-o(1)}
 \le |A|\le N^{1/2+(1/2-\theta)+o(1)}.
\]
The gap between the upper and lower exponents is \(1-2\theta>0\). As \(\theta\) approaches \(1/2\) from below, both exponents approach \(1/2\).
If we subtract the upper exponent from the lower one and leave out the \(o(1)\) terms, we get
\[
 \theta-(1-\theta)=2\theta-1=2\left(\theta-\frac12\right).
\]
The lower exponent is larger than the upper one exactly when \(\theta>1/2\). Thus these bounds give no contradiction for \(\theta<1/2\). We need stronger bounds when \(\theta\) is close to \(1/2\) from below. Suppose that we can prove
\[
 N^{k_1\theta-k_2-o(1)}\le |A|\le N^{k_3-k_4\theta+o(1)},
\]
where \(k_1,k_2,k_3,k_4\) are fixed real constants and \(k_1+k_4>0\). Subtracting the upper exponent from the lower one now gives
\[
 \begin{aligned}
 \Delta_\theta
 &=(k_1\theta-k_2-o(1))-(k_3-k_4\theta+o(1))\\
 &=(k_1+k_4)\theta-(k_2+k_3)-o(1).
 \end{aligned}
\]
If
\[
 \theta>\frac12-\left(\frac12-\frac{k_2+k_3}{k_1+k_4}\right)
       =\frac{k_2+k_3}{k_1+k_4},
\]
then \(\Delta_\theta>0\) for all sufficiently large \(N\). The lower bound for \(|A|\) would then exceed the upper bound. To make this happen for some \(\theta<1/2\), we look for \((k_1,k_2,k_3,k_4)\) such that
\[
 k_1+k_4>2(k_2+k_3).
\]
\end{remark}

\begin{remark}
We now turn to chords of the cubic curve. For \(a,b\in A\) with \(a>b\), consider the chord from \((b,b^3)\) to \((a,a^3)\). If a positive integer \(p\) divides \(a-b\), its chord vector divided by \(p\) is an integer point:
\[
 \left(\frac{a-b}{p},\,\frac{a^3-b^3}{p}\right)\in\mathbb Z^2.
\]
Propositions \ref{th:chord-ratio} and \ref{th:chords} show that a large eighth moment produces a large set of such points in a short horizontal interval, with a small difference set.
\end{remark}

\begin{proposition}\label{th:chord-ratio}
Let \(13/27<\theta<1/2\), and let \(N\) be sufficiently large.
If a nonempty set \(A\subseteq\{1,\ldots,N\}\) satisfies \(I_8(\mathbf1_A)\ge N^\theta|A|^4\), then there exist a prime \(p\), a scale \(X=2^j\) with \(j\) a nonnegative integer, and a nonempty set
\[
 E\subseteq\left\{\left(\frac{a-b}{p},\frac{a^3-b^3}{p}\right):
 \begin{array}{l}
 a,b\in A,\ a>b,\ p\mid(a-b),\\
 X\le(a-b)/p<2X
 \end{array}\right\}
\]
such that
\[
 |E|\ge N^{2\theta-1/2-o(1)}.
\]
For every integer \(1\le H\le X\), partition \([X,2X)\) into
\[
 I_j=[X+jH,\,\min\{X+(j+1)H,\,2X\}),\qquad
 j=0,\ldots,\lceil X/H\rceil-1.
\]
If \(h_0\) is the left endpoint of an interval maximizing
\(\#\{(h,k)\in E:h\in I_j\}\), then the nonempty set
\[
 B=\{(h,k)\in E:h_0\le h<h_0+H\}
\]
satisfies
\[
 \frac{|B-B|}{|B|}\le N^{153(1/2-\theta)+o(1)}.
\]
The choices of \(p,X,E\) are independent of \(H\), and the bound is uniform in \(H\) and in the choice of a maximizing interval.
\end{proposition}

We next choose \(H\) and obtain a lower bound for \(|B|\). The estimate for \(|B-B|/|B|\) still holds for this same \(B\).

\begin{proposition}\label{th:chords}
Let \(1/2-1/1550<\theta<1/2\) and let \(N\) be sufficiently large. Suppose that a nonempty set \(A\subseteq\{1,\ldots,N\}\) satisfies \(I_8(\mathbf1_A)\ge N^\theta|A|^4\). Let \(p,X,E\) be as in Proposition~\ref{th:chord-ratio}. Then
\[
 X\ge N^{1/2-155(1/2-\theta)-o(1)}.
\]
For sufficiently large \(N\), the integer
\[
 H=\lfloor XN^{-2/5}\rfloor
\]
satisfies \(1\le H\le XN^{-2/5}\). At this width, the set \(B\) selected in Proposition~\ref{th:chord-ratio} satisfies
\[
 |B|\ge N^{2\theta-9/10-o(1)}.
\]
Its horizontal coordinates lie in \([h_0,h_0+H]\), with \(h_0\ge X\), and
\[
 B\subseteq\left\{\left(\frac{a-b}{p},\frac{a^3-b^3}{p}\right):
 a,b\in A,\ a>b,\ p\mid(a-b)\right\}.
\]
For this same set \(B\), the difference-set estimate in Proposition~\ref{th:chord-ratio} remains valid.
\end{proposition}

\begin{remark}
In Propositions \ref{th:chord-ratio} and \ref{th:chords}, as \(\theta\) approaches \(1/2\) from below,
\[
 2\theta-\frac9{10}\longrightarrow\frac1{10},
 \qquad 153\left(\frac12-\theta\right)\longrightarrow0.
\]
Thus the lower bound for \(|B|\) retains a fixed positive power of \(N\), while the power in the upper bound for \(|B-B|/|B|\) approaches zero.
The next proposition bounds \(|B-B|/|B|\) from below in terms of \(|B|\). Comparing this bound with Propositions~\ref{th:chord-ratio} and \ref{th:chords} will give the contradiction.
\end{remark}

\begin{proposition}\label{th:growth}
Let \(N\ge2\) and \(p\ge1\) be integers, and let \(B\) be any nonempty subset of the following set of integer chord points,
\[
 B\subseteq\left\{\left(\frac{a-b}{p},\,\frac{a^3-b^3}{p}\right):
       a,b\in\mathbb Z,\ 1\le b<a\le N,\ p\mid(a-b)\right\}.
\]
Let \(1\le H\le X\le h_0\). If all horizontal coordinates of \(B\) lie in \([h_0,h_0+H]\) and \(H\le XN^{-2/5}\), then, for every \(\eta>0\), there is a constant \(C_\eta>0\), depending only on \(\eta\), such that
\[
 \frac{|B-B|}{|B|}\ge C_\eta N^{-\eta}|B|^{1/31}.
\]
\end{proposition}

We can now prove Theorem~\ref{th:main}.

\begin{proof}[Proof of Theorem \ref{th:main}]
Fix
\[
 \frac12-\frac1{47500}<\theta<\frac12 .
\]
We first prove the estimate for indicator functions. Suppose, to the contrary, that for arbitrarily large \(N\) there is a nonempty \(A\subseteq\{1,\ldots,N\}\) with \(I_8(\mathbf1_A)\ge N^\theta|A|^4\). Propositions~\ref{th:chord-ratio} and~\ref{th:chords} give a nonempty set \(B\) satisfying
\[
 |B|\ge N^{2\theta-9/10-o(1)},\qquad
 \frac{|B-B|}{|B|}\le N^{153(1/2-\theta)+o(1)},
\]
\[
 1\le H=\lfloor XN^{-2/5}\rfloor\le XN^{-2/5}.
\]
The width condition allows us to apply Proposition~\ref{th:growth} to this same \(B\). For every fixed \(\eta>0\), the lower and upper bounds for its difference-set ratio give
\[
 C_\eta N^{-\eta+(2\theta-9/10-o(1))/31}
 \le\frac{|B-B|}{|B|}
 \le N^{153(1/2-\theta)+o(1)} .
\]
Raise the two outer expressions to the \(31\)-st power and bring the powers of \(N\) to the left. Since \(31\cdot153=4743\),
\begin{equation}\label{chord:eq:contradiction}
 C_\eta^{31}
 N^{2\theta-9/10-4743(1/2-\theta)-31\eta-o(1)}\le1 .
\end{equation}
Before subtracting \(31\eta+o(1)\), the exponent in \eqref{chord:eq:contradiction} is positive whenever \(1/2-\theta<1/(10\cdot4745)=1/47450\). The saving \(1/47500\) in Theorem~\ref{th:main} leaves the following margin.
\[
 \begin{aligned}
 2\theta-\frac9{10}-4743\left(\frac12-\theta\right)
 &=\frac1{10}-4745\left(\frac12-\theta\right)\\
 &>\frac1{10}-\frac{4745}{47500}
 =\frac1{9500}>0 .
 \end{aligned}
\]
Choose \(\eta>0\) small enough that
\[
 0<31\eta<
 \frac12\left(\frac1{10}-4745\left(\frac12-\theta\right)\right).
\]
With this fixed choice of \(\eta\), the exponent in \eqref{chord:eq:contradiction} stays positive as \(N\to\infty\). The left-hand side therefore tends to infinity, contradicting \eqref{chord:eq:contradiction}. Hence, for all sufficiently large \(N\), every nonempty \(A\subseteq\{1,\ldots,N\}\) satisfies \(I_8(\mathbf1_A)<N^\theta|A|^4\).
Only finitely many \(N\) remain. For them, use \(I_8(\mathbf1_A)\le|A|^8\) and increase the constant. We obtain
\[
 I_8(\mathbf1_A)\le C_\theta N^\theta|A|^4
 \qquad(N\ge1,\ A\subseteq\{1,\ldots,N\}).
\]

Given \(\varepsilon>0\), choose
\[
 \frac12-\frac1{47500}<\theta<
 \min\left\{\frac12,\frac12-\frac1{47500}+\frac{\varepsilon}{2}\right\}.
\]
By the indicator estimate \eqref{eq:indicator-target} just proved, we may take \(D=C_\theta N^\theta\) in \eqref{eq:weighted-reduction}. The bound \((1+\log N)^4\ll_\varepsilon N^{\varepsilon/2}\) absorbs the logarithmic factor into the remaining power of \(N\), and gives
\[
 I_8(b)\le C_\theta N^\theta(1+\log N)^4
       \left(\sum_{n=1}^N|b_n|^2\right)^4
 \le C_\varepsilon N^{1/2-1/47500+\varepsilon}
       \left(\sum_{n=1}^N|b_n|^2\right)^4 .
\]
This completes the proof of Theorem~\ref{th:main}.
\end{proof}

We prove Proposition \ref{th:support} in Section~\ref{sec:support}, and Propositions \ref{th:chord-ratio} and \ref{th:chords} in Section~\ref{sec:chords}.
In Section~\ref{sec:growth}, we first prove a difference-set bound for points with integer slopes and then apply it to chord points on monomial curves. Taking the curve degree \(d=3\) and the difference order \(r=5\) gives Proposition \ref{th:growth}.

\section{Proof of Proposition \ref*{th:support}}\label{sec:support}
\begin{lemma}\label{le:vertical-cubic}
Let \(N\ge2\) and \(p\ge1\) be integers, and let
\[
 U\subseteq\left\{\left(\frac{a-b}{p},\frac{a^3-b^3}{p}\right):
 a,b\in\mathbb Z,\ 1\le b<a\le N,\ p\mid(a-b)\right\}.
\]
Every point \((h,k)\in U\) uniquely determines its endpoints \(a>b\). For every \(t\in\mathbb Z\setminus\{0\}\),
\[
 \#\{(h,k)\in U:(h,k+t)\in U\}\le 2\tau(|t|)^2.
\]
Moreover,
\[
 \#\{(h,k)\in U:(h,k+t)\in U\}\le N^{o(1)},
\]
uniformly in \(p,t,U\).
\end{lemma}

\begin{proof}
Take integer endpoints with \(1\le b<a\le N\) and \(p\mid(a-b)\). Their chord point satisfies
\[
 h=\frac{a-b}{p},\qquad
 k=\frac{a^3-b^3}{p}=\frac{(b+ph)^3-b^3}{p}
   =3hb^2+3ph^2b+p^2h^3 .
\]
For fixed \(p,h>0\), \(3hb^2+3ph^2b+p^2h^3\) is strictly increasing in \(b\ge1\).
Once \(p,h,k\) are fixed, there is therefore at most one positive integer \(b\). The value of \(a\) is then fixed by \(a=b+ph\).

Let
\[
 U\subseteq\left\{\left(\frac{a-b}{p},\frac{a^3-b^3}{p}\right):
 a,b\in\mathbb Z,\ 1\le b<a\le N,\ p\mid(a-b)\right\}.
\]
Fix \(t\in\mathbb Z\setminus\{0\}\). For \((h,k),(h,k+t)\in U\), write the smaller endpoints as \(b\) and \(b+w\), respectively. Since both chord points have horizontal coordinate \(h\),
\[
 \begin{aligned}
 (h,k)&=\left(h,\frac{(b+ph)^3-b^3}{p}\right),\\
 (h,k+t)&=\left(h,\frac{(b+w+ph)^3-(b+w)^3}{p}\right).
 \end{aligned}
\]
Subtracting their vertical coordinates gives
\[
 \begin{aligned}
 t&=\frac{(b+w+ph)^3-(b+w)^3}{p}-\frac{(b+ph)^3-b^3}{p}\\
  &=3h\bigl((b+w)^2-b^2\bigr)+3ph^2w
   =3hw(2b+w+ph).
 \end{aligned}
\]
Since \(t\ne0\),
\[
 h\mid |t|,\qquad h>0,\qquad w\ne0,\qquad |w|\mid |t| .
\]
For fixed \(h,w\), this equation determines \(b\), and then \(k\):
\[
 b=\tfrac12\left(\frac{t}{3hw}-w-ph\right),\qquad
 k=\frac{(b+ph)^3-b^3}{p}.
\]
Thus each choice of \((h,w)\) gives at most one point \((h,k)\). There are \(\tau(|t|)\) possible positive divisors \(h\), and at most \(2\tau(|t|)\) possible signed divisors \(w\). Some choices may give no admissible endpoints, but counting all of them gives the upper bound
\begin{equation}\label{eq:vertical}
 \begin{aligned}
 \#\{(h,k)\in U:(h,k+t)\in U\}
 &\le \#\{h>0:h\mid|t|\}\,
     \#\{w\ne0:|w|\mid|t|\}\\
 &=2\tau(|t|)^2 .
 \end{aligned}
\end{equation}
If a pair occurs, its vertical difference has absolute value at most \(N^3\). The divisor bound in \eqref{eq:vertical} therefore gives the uniform estimate; if \(|t|>N^3\), the count is zero:
\[
 \begin{gathered}
 (h,k),(h,k+t)\in U
 \quad\Longrightarrow\quad
 0<k,\ k+t\le N^3
 \quad\Longrightarrow\quad 0<|t|\le N^3,\\
 0<|t|\le N^3\quad\Longrightarrow\quad
 2\tau(|t|)^2\le N^{o(1)},\\
 |t|>N^3\quad\Longrightarrow\quad
 \#\{(h,k)\in U:(h,k+t)\in U\}=0 .
 \end{gathered}
\]

\end{proof}

Proposition~\ref{chord:le:vertical} extends the vertical-difference bound in Lemma~\ref{le:vertical-cubic} to integer polynomials of degree at least three.

\begin{restateproposition}{th:support}
Let \(13/27<\theta<1/2\), and let \(A\) be a nonempty subset of \(\{1,\ldots,N\}\).
If \(I_8(\mathbf1_A)\ge N^\theta|A|^4\), then
\[
 N^{\theta-o(1)}\le |A|\le N^{1-\theta+o(1)} .
\]
\end{restateproposition}

\begin{proof}[Proof of Proposition \ref{th:support}]
We begin with the number of triples representing a fixed \((s,v)\in\mathbb Z^2\). If \(a+b+c=s,\ a^3+b^3+c^3=v\) has no solution in \(A^3\), the count is zero. Otherwise, take any solution \((a,b,c)\). Substituting into \(s^3-v\), and expanding \(((a+b)+c)^3\) followed by \((a+b)^3\), gives
\[
 \begin{aligned}
 s^3-v
 &=(a+b+c)^3-a^3-b^3-c^3\\
 &=3c(a+b)(a+b+c)+3ab(a+b)\\
 &=3(a+b)\bigl(ab+ac+bc+c^2\bigr)\\
 &=3(a+b)(a+c)(b+c).
 \end{aligned}
\]
Since \(1\le a,b,c\le N\),
\[
 \begin{gathered}
 2\le a+b,\ a+c,\ b+c\le2N,\\
 8\le \frac{s^3-v}{3}=(a+b)(a+c)(b+c)\le (2N)^3=8N^3 .
 \end{gathered}
\]
Thus \((s^3-v)/3\) is a positive integer fixed by \((s,v)\). The product \((a+b)(a+c)(b+c)\) is fixed, so choosing any two of \(a+b,\ a+c,\ b+c\) determines the third. For example,
\[
 b+c=\frac{s^3-v}{3(a+b)(a+c)} .
\]
Each of \(a+b\) and \(a+c\) must be a positive divisor of \((s^3-v)/3\). Hence there are at most \(\tau((s^3-v)/3)^2\) choices for the ordered triple \((a+b,a+c,b+c)\). Such a choice determines at most one \((a,b,c)\), since
\[
 \begin{aligned}
 a&=\tfrac12\bigl((a+b)+(a+c)-(b+c)\bigr),\\
 b&=\tfrac12\bigl((a+b)+(b+c)-(a+c)\bigr),\\
 c&=\tfrac12\bigl((a+c)+(b+c)-(a+b)\bigr).
 \end{aligned}
\]
Fix \(\eta>0\) and use the divisor bound with exponent \(\eta/6\). Squaring doubles this exponent, and \((s^3-v)/3\le8N^3\), so
\(\tau((s^3-v)/3)^2\ll_\eta (8N^3)^{\eta/3}\ll_\eta N^\eta\). The implied constant is independent of \((s,v)\). Since \(\eta>0\) is arbitrary, and a pair \((s,v)\) with no representation contributes zero, we obtain
\begin{equation}\label{eq:triples}
 \#\{(a,b,c)\in A^3:a+b+c=s,\ a^3+b^3+c^3=v\}\le N^{o(1)} .
\end{equation}

We now use this count to bound \(I_6(\mathbf1_A)\). Write the sixth power of the modulus as the cube of the sum times the cube of its complex conjugate:
\[
 \left|\sum_{a\in A}e^{2\pi i(ax+a^3y)}\right|^6=
 \left(\sum_{a\in A}e^{2\pi i(ax+a^3y)}\right)^3
 \left(\sum_{c\in A}e^{-2\pi i(cx+c^3y)}\right)^3 .
\]
Conjugation changes the sign in the exponential. Expanding the two cubes introduces three variables \(a_1,a_2,a_3\) and three variables \(c_1,c_2,c_3\). Integrating the resulting finite sum term by term gives
\[
 \begin{aligned}
 I_6(\mathbf1_A)
 &=\sum_{\substack{a_1,a_2,a_3\in A\\c_1,c_2,c_3\in A}}
 \left(\int_0^1
 \exp\!\bigl(2\pi i x(a_1+a_2+a_3-c_1-c_2-c_3)\bigr)\,dx\right)\\
 &\qquad{}\times\left(\int_0^1
 \exp\!\bigl(2\pi i y(a_1^3+a_2^3+a_3^3-c_1^3-c_2^3-c_3^3)\bigr)\,dy\right).
 \end{aligned}
\]
For every integer \(m\),
\begin{equation}\label{eq:orthogonality}
 \int_0^1 e^{2\pi i m x}\,dx=
 \begin{cases}
 1,&m=0,\\
 (e^{2\pi i m}-1)/(2\pi i m)=0,&m\ne0.
 \end{cases}
\end{equation}
By \eqref{eq:orthogonality}, a sextuple contributes \(1\) exactly when the two triples have the same sum and the same sum of cubes. Writing these common values as \(s,v\), we have
\[
 \begin{aligned}
 s&=a_1+a_2+a_3=c_1+c_2+c_3,\\
 v&=a_1^3+a_2^3+a_3^3=c_1^3+c_2^3+c_3^3.
 \end{aligned}
\]
Every other sextuple contributes \(0\). Grouping the surviving sextuples by \((s,v)\), we obtain
\begin{equation}\label{eq:six-solutions}
 I_6(\mathbf1_A)=
 \sum_{(s,v)\in\mathbb Z^2}
 \bigl(\#\{(a,b,c)\in A^3:a+b+c=s,\ a^3+b^3+c^3=v\}\bigr)^2 .
\end{equation}
By \eqref{eq:triples},
\[
 I_6(\mathbf1_A)\le N^{o(1)}
 \sum_{(s,v)\in\mathbb Z^2}
 \#\{(a,b,c)\in A^3:a+b+c=s,\ a^3+b^3+c^3=v\}.
\]
Each \((a,b,c)\in A^3\) determines exactly one pair \((s,v)\). Counting first by the triple instead of by \((s,v)\), we get
\[
 \begin{aligned}
 &\sum_{(s,v)\in\mathbb Z^2}
 \#\{(a,b,c)\in A^3:a+b+c=s,\ a^3+b^3+c^3=v\}\\
 &=\sum_{(a,b,c)\in A^3}
 \#\{(s,v)\in\mathbb Z^2:(s,v)=(a+b+c,a^3+b^3+c^3)\}\\
 &=\sum_{a\in A}\sum_{b\in A}\sum_{c\in A}1=|A|^3 .
 \end{aligned}
\]
For every \((x,y)\in[0,1]^2\), the triangle inequality gives
\[
 \left|\sum_{a\in A} e^{2\pi i(ax+a^3y)}\right|
 \le\sum_{a\in A}\left|e^{2\pi i(ax+a^3y)}\right|
 =\sum_{a\in A}1=|A|.
\]
Consequently,
\begin{equation}\label{eq:six-eight}
 \begin{aligned}
 I_6(\mathbf1_A)&\le N^{o(1)}|A|^3,\\
 I_8(\mathbf1_A)
 &=\int_{[0,1]^2}
 \left|\sum_{a\in A}e^{2\pi i(ax+a^3y)}\right|^2
 \left|\sum_{a\in A}e^{2\pi i(ax+a^3y)}\right|^6\,dx\,dy\\
 &\le |A|^2 I_6(\mathbf1_A)\le N^{o(1)}|A|^5 .
 \end{aligned}
\end{equation}
By \eqref{eq:six-eight} and the hypothesis of Proposition~\ref{th:support},
\[
 N^\theta|A|^4\le I_8(\mathbf1_A)\le N^{o(1)}|A|^5,\qquad |A|>0
 \quad\Longrightarrow\quad |A|\ge N^{\theta-o(1)} .
\]

We next prove the upper bound for \(|A|\). We construct a graph on the triple sums and use Lemma~\ref{le:vertical-cubic} to bound the number of common neighbors of two vertices in the same vertical column.

Use the set
\[
 \mathcal T=\{(a+b+c,a^3+b^3+c^3):a,b,c\in A\}
\]
as the vertex set. Join two vertices when their difference is a nonzero chord vector \((a-b,a^3-b^3)\), with \(a,b\in A\). Swapping \(a,b\) reverses the vector, so adjacency is symmetric. In counting edges, we count both orientations separately.

Fix \(P=(s,v)\) and \(P'=(s,v')\) with \(v\ne v'\). Their equal horizontal coordinates will turn the common-neighbor problem into a vertical-difference count. Suppose that \(Q=(Q_1,Q_2)\) is adjacent to both. We must have \(Q_1\ne s\), since a chord vector with zero horizontal coordinate is zero. Choose the orientation of both chord vectors so that their horizontal coordinates are positive. According to which side of the column contains \(Q\), we obtain
\[
 \begin{array}{ll}
 Q_1<s:&(P'-Q)-(P-Q)=(0,v'-v),\\
 Q_1>s:&(Q-P')-(Q-P)=(0,v-v').
 \end{array}
\]
When \(Q_1<s\), apply \eqref{eq:vertical} to \(P-Q,P'-Q\), whose vertical difference is \(v'-v\). When \(Q_1>s\), apply it to \(Q-P,Q-P'\), whose vertical difference is \(v-v'\). In each case \(p=1\), and the chord vector involving \(P\) determines \(Q\) uniquely. Hence
\begin{equation}\label{eq:common-neighbors}
 \begin{aligned}
 &\#\{Q\in\mathcal T:Q\text{ is adjacent to both }P,P'\}\\
 &=\#\{Q\in\mathcal T:Q_1<s,\ Q\text{ is adjacent to both }P,P'\}\\
 &\quad+\#\{Q\in\mathcal T:Q_1>s,\ Q\text{ is adjacent to both }P,P'\}\\
 &\le 2\tau(|v'-v|)^2+2\tau(|v-v'|)^2
 \le N^{o(1)} .
 \end{aligned}
\end{equation}
To justify the final bound, note that \(P,P'\in\mathcal T\) and \(A\subseteq\{1,\ldots,N\}\) imply
\[
 3\le v,v'\le3N^3,\qquad v\ne v'
 \quad\Longrightarrow\quad
 1\le |v'-v|=|v-v'|\le3N^3 .
\]
For every \(\eta>0\), taking \(\varepsilon=\eta/6\) in \(\tau(m)\ll_\varepsilon m^\varepsilon\) gives
\[
 \tau(m)\le C_{\eta/6}m^{\eta/6}
 \qquad(m\in\mathbb Z_{\ge1}).
\]
Therefore
\[
 \begin{aligned}
 2\tau(|v'-v|)^2+2\tau(|v-v'|)^2
 &=4\tau(|v'-v|)^2\\
 &\le4C_{\eta/6}^{\,2}|v'-v|^{\eta/3}\\
 &\le4C_{\eta/6}^{\,2}(3N^3)^{\eta/3}\\
 &=4\cdot3^{\eta/3}C_{\eta/6}^{\,2}N^\eta
 \ll_\eta N^\eta .
 \end{aligned}
\]
This holds for every \(\eta>0\), with a constant independent of \(P,P'\). We may therefore write the common-neighbor bound as \(N^{o(1)}\).

We now use \eqref{eq:common-neighbors} to bound the number of edges. Write
\[
 \mathcal T_s=\{(s',v)\in\mathcal T:s'=s\},\qquad
 \deg_s(Q)=\#\{P\in\mathcal T_s:P\text{ is adjacent to }Q\}.
\]
Count the ordered triples \((P,P',Q)\) with distinct \(P,P'\in\mathcal T_s\) and with \(Q\) adjacent to both. Counting first by \(Q\) and then by \((P,P')\) gives
\[
 \begin{aligned}
 \sum_{Q\in\mathcal T}\deg_s(Q)(\deg_s(Q)-1)
 &=\sum_{\substack{P,P'\in\mathcal T_s\\P\ne P'}}
   \#\{Q\in\mathcal T:Q\text{ is adjacent to both }P,P'\}\\
 &\le N^{o(1)}|\mathcal T_s|(|\mathcal T_s|-1).
 \end{aligned}
\]
Using
\[
 j\le1+\sqrt{j(j-1)}\qquad(j\in\mathbb Z_{\ge0})
\]
and applying the Cauchy--Schwarz inequality followed by \eqref{eq:common-neighbors}, we obtain
\[
 \begin{aligned}
 \sum_{Q\in\mathcal T}\deg_s(Q)
 &\le |\mathcal T|+
   |\mathcal T|^{1/2}
   \left(\sum_{Q\in\mathcal T}\deg_s(Q)(\deg_s(Q)-1)\right)^{1/2}\\
 &\le |\mathcal T|+N^{o(1)}|\mathcal T|^{1/2}|\mathcal T_s|.
 \end{aligned}
\]
Since
\[
 \#\{s:\mathcal T_s\ne\varnothing\}\le3N,\qquad
 \sum_s|\mathcal T_s|=|\mathcal T|\le|A|^3,
\]
the sum over nonempty columns gives the following bound for the number of directed edges:
\begin{equation}\label{eq:edge-bound}
 \begin{aligned}
 \sum_{s:\,\mathcal T_s\ne\varnothing}\ \sum_{Q\in\mathcal T}\deg_s(Q)
 &\le 3N|\mathcal T|+N^{o(1)}|\mathcal T|^{3/2}\\
 &\le 3N|A|^3+N^{o(1)}|A|^{9/2}.
 \end{aligned}
\end{equation}

We now relate these edges to \(I_8(\mathbf1_A)\). Expanding four copies of the exponential sum and four copies of its complex conjugate, and using \eqref{eq:orthogonality}, gives the identity \eqref{eq:eight-solutions} in the form
\[
 I_8(\mathbf1_A)=
 \#\left\{(a_1,\ldots,a_4,c_1,\ldots,c_4)\in A^8:
 \begin{array}{l}
 \sum_{j=1}^{4}(a_j-c_j)=0,\\
 \sum_{j=1}^{4}(a_j^3-c_j^3)=0
 \end{array}\right\}.
\]
First consider \(a_4=c_4\). The terms \(a_4,c_4\) cancel in the linear equation, and \(a_4^3,c_4^3\) cancel in the cubic equation. There are \(|A|\) choices for their common value. The remaining six variables are counted by \(I_6(\mathbf1_A)\) by \eqref{eq:six-solutions}, so this case contributes \(|A|I_6(\mathbf1_A)\) solutions.

Now suppose that \(a_4\ne c_4\), and take
\[
 \begin{aligned}
 P&=(a_1+a_2+a_3,\,a_1^3+a_2^3+a_3^3),\\
 Q&=(c_1+c_2+c_3,\,c_1^3+c_2^3+c_3^3).
 \end{aligned}
\]
Since \(\sum_{j=1}^4(a_j-c_j)=\sum_{j=1}^4(a_j^3-c_j^3)=0\),
\[
 P-Q=(c_4-a_4,\ c_4^3-a_4^3).
\]
As \(P,Q\in\mathcal T\), \(a_4,c_4\in A\), and \(a_4\ne c_4\),
\[
 P-Q=(c_4-a_4,c_4^3-a_4^3)
 \in\{(a-b,a^3-b^3):a,b\in A,\ a\ne b\}.
\]
Hence \(P\) and \(Q\) are adjacent, and \((P,Q)\) is one of the directed edges we have counted. The endpoints of the corresponding chord are
\[
 \begin{array}{ll}
 c_4>a_4:&P-Q=(c_4,c_4^3)-(a_4,a_4^3),\\[3pt]
 a_4>c_4:&Q-P=(a_4,a_4^3)-(c_4,c_4^3).
 \end{array}
\]
In either case, the positively oriented chord vector determines the two endpoints by Lemma~\ref{le:vertical-cubic} with \(p=1\). The sign of \(P_1-Q_1\) fixes the order of \(a_4,c_4\).

Thus fixing the directed edge \((P,Q)\) leaves at most one choice for \((a_4,c_4)\):
\[
 \#\left\{(a_4,c_4)\in A^2:
 \begin{array}{l}
 c_4-a_4=P_1-Q_1,\\
 c_4^3-a_4^3=P_2-Q_2
 \end{array}\right\}\le1.
\]
For fixed \(P,Q\), each of the two ordered triples representing them has at most \(N^{o(1)}\) choices by \eqref{eq:triples}. Together they have at most \(N^{o(1)}\cdot N^{o(1)}=N^{o(1)}\) choices. Adding the contribution with \(a_4=c_4\), and using \eqref{eq:six-eight} and \eqref{eq:edge-bound}, gives
\[
 \begin{aligned}
 I_8(\mathbf1_A)
 &\le |A|I_6(\mathbf1_A)
      +N^{o(1)}\bigl(3N|A|^3+N^{o(1)}|A|^{9/2}\bigr)\\
 &\le N^{o(1)}\bigl(|A|^4+N|A|^3+|A|^{9/2}\bigr).
 \end{aligned}
\]
Dividing by \(|A|^4\) gives
\begin{equation}\label{eq:graph-moment}
 \frac{I_8(\mathbf1_A)}{|A|^4}
 \le N^{o(1)}\left(1+\frac{N}{|A|}+|A|^{1/2}\right).
\end{equation}

The term \(|A|^{1/2}\) in \eqref{eq:graph-moment} can be too large when \(|A|\) is large. In that range we use Wooley's estimate on the full interval \cite[Theorem 1.1]{wooley}:
\begin{equation}\label{eq:wooley-interval}
 I_8(\mathbf1_{\{1,\ldots,N\}})\le N^{13/3+o(1)}.
\end{equation}
Every solution with its eight entries in \(A\) is also counted on \(\{1,\ldots,N\}\). The bounds \(|A|^{1/2}\) and \(N^{13/3}/|A|^4\) agree when \(|A|=N^{26/27}\). We therefore split at this value. If \(|A|\le N^{26/27}\), then \(|A|^{1/2}\le N^{13/27}\). If \(|A|\ge N^{26/27}\), \eqref{eq:wooley-interval} gives
\[
 \frac{I_8(\mathbf1_A)}{|A|^4}
 \le\frac{I_8(\mathbf1_{\{1,\ldots,N\}})}{|A|^4}
 \le N^{13/3-104/27+o(1)}
 =N^{13/27+o(1)} .
\]
Combining \eqref{eq:graph-moment} for \(|A|\le N^{26/27}\) with the bound from \eqref{eq:wooley-interval} for \(|A|\ge N^{26/27}\), and using \(I_8(\mathbf1_A)\ge N^\theta|A|^4\), we obtain
\[
 N^\theta
 \le\frac{I_8(\mathbf1_A)}{|A|^4}
 \le N^{o(1)}\left(N^{13/27}+\frac N{|A|}\right).
\]
Since \(\theta>13/27\) is fixed,
\[
 \frac{N^{13/27+o(1)}}{N^\theta}
 =N^{-(\theta-13/27)+o(1)}\longrightarrow0 .
\]
We may therefore absorb the \(N^{13/27+o(1)}\) term into half of \(N^\theta\). For sufficiently large \(N\),
\[
 \frac{N^\theta}{2}\le\frac{N^{1+o(1)}}{|A|}
 \quad\Longrightarrow\quad
 |A|\le N^{1-\theta+o(1)} .
\]
Together with the lower bound \(|A|\ge N^{\theta-o(1)}\), this proves the assertion.
\end{proof}

\begin{remark}
Proposition \ref{th:support} shows that an eighth moment of size at least \(N^\theta|A|^4\) forces both
\[
 |A|\ge N^{\theta-o(1)},\qquad |A|\le N^{1-\theta+o(1)}.
\]
The lower bound comes from counting triples in the sixth moment. For the upper bound, we count common neighbors in the graph and then use the eighth-moment estimate on the full interval.
Thus, when \(\theta\) is close to \(1/2\), the remaining argument only has to treat \(|A|\) between powers of \(N\) whose exponents are close to \(1/2\).
\end{remark}

\newpage
\section{Proofs of Propositions \ref*{th:chord-ratio} and \ref*{th:chords}}\label{sec:chords}
\begin{lemma}\label{le:moment-extraction}
Let \(13/27<\theta<1/2\), and suppose that a nonempty set \(A\subseteq\{1,\ldots,N\}\) satisfies \(I_8(\mathbf1_A)\ge N^\theta|A|^4\). There exist nonempty sets
\[
 \Omega\subseteq\{(a+b+c,a^3+b^3+c^3):a,b,c\in A\},
 \qquad S\subseteq A,
\]
and \((u_0,v_0)\in A^2\) such that
\[
\begin{aligned}
 |\Omega|&\ge N^{3\theta/2-o(1)}|A|^{3/2},&
 \frac{|3\Omega-3\Omega|}{|\Omega|}&\le N^{144(1/2-\theta)+o(1)},\\
 |S|&\ge N^{2\theta-1/2-o(1)},&
 \frac{|A|}{|S|}&\le N^{3(1/2-\theta)+o(1)} .
\end{aligned}
\]
For \(a\in A\), call \((u,v)\in A^2\) a partner pair of \(a\) if \((a+u+v,a^3+u^3+v^3)\in\Omega\). Each \(a\in S\) has \((u_0,v_0)\) as a partner pair and has at least \(|\Omega|/(2|A|)\) partner pairs.
\end{lemma}

\begin{proof}
The hypothesis \(13/27<\theta<1/2\) allows us to apply Proposition~\ref{th:support}. Use the set \(\mathcal T=\{(a+b+c,a^3+b^3+c^3):a,b,c\in A\}\). We first show that it contains many additive quadruples. H\"older's inequality gives
\[
 \begin{aligned}
 I_8(\mathbf1_A)
 &=\int_{[0,1]^2}
 \left|\sum_{a\in A}e^{2\pi i(ax+a^3y)}\right|^4
 \left|\sum_{a\in A}e^{2\pi i(ax+a^3y)}\right|^4\,dx\,dy\\
 &\le I_6(\mathbf1_A)^{2/3}I_{12}(\mathbf1_A)^{1/3}.
 \end{aligned}
\]
By \eqref{eq:six-solutions},
\[
 I_6(\mathbf1_A)\ge
 \#\{(a_1,a_2,a_3,c_1,c_2,c_3)\in A^6:a_j=c_j,\ 1\le j\le3\}
 =|A|^3>0 .
\]
The diagonal solutions show that \(I_6(\mathbf1_A)>0\), so we can rearrange the preceding moment inequality to obtain
\begin{equation}\label{eq:twelve-lower}
 I_{12}(\mathbf1_A)\ge
 \frac{I_8(\mathbf1_A)^3}{I_6(\mathbf1_A)^2}.
\end{equation}

Applying \eqref{eq:orthogonality} to the twelfth moment shows that \(I_{12}(\mathbf1_A)\) counts the ordered twelve-tuples satisfying \(\sum_{j=1}^6a_j=\sum_{j=1}^6c_j,\) \(\sum_{j=1}^6a_j^3=\sum_{j=1}^6c_j^3\). Group the variables into four triples and associate to them the points
\[
 \begin{aligned}
 P_1&=(a_1+a_2+a_3,\,a_1^3+a_2^3+a_3^3),&
 P_2&=(a_4+a_5+a_6,\,a_4^3+a_5^3+a_6^3),\\
 P_3&=(c_1+c_2+c_3,\,c_1^3+c_2^3+c_3^3),&
 P_4&=(c_4+c_5+c_6,\,c_4^3+c_5^3+c_6^3).
 \end{aligned}
\]
\[
 \left.
 \begin{aligned}
 \sum_{j=1}^6a_j&=\sum_{j=1}^6c_j\\
 \sum_{j=1}^6a_j^3&=\sum_{j=1}^6c_j^3
 \end{aligned}
 \right\}
 \quad\Longleftrightarrow\quad P_1+P_2=P_3+P_4.
\]
For each additive quadruple, the number of twelve-tuples representing it is the product of the four triple representation counts. Thus
\[
 I_{12}(\mathbf1_A)=
 \sum_{\substack{P_1,P_2,P_3,P_4\in\mathcal T\\P_1+P_2=P_3+P_4}}
 \ \prod_{j=1}^4\#\{(a,b,c)\in A^3:(a+b+c,a^3+b^3+c^3)=P_j\}.
\]
By \eqref{eq:triples}, this product is at most \((N^{o(1)})^4=N^{o(1)}\). Using \eqref{eq:twelve-lower}, the sixth-moment bound in \eqref{eq:six-eight}, and the assumption \(I_8(\mathbf1_A)\ge N^\theta|A|^4\), we obtain
\begin{equation}\label{eq:energy}
 \begin{aligned}
 &\#\{(P_1,P_2,P_3,P_4)\in\mathcal T^4:P_1+P_2=P_3+P_4\}\\
 &\qquad\ge N^{-o(1)} I_{12}(\mathbf1_A)
 \ge N^{-o(1)}\frac{I_8(\mathbf1_A)^3}{|A|^6}
 \ge N^{3\theta-o(1)}|A|^6 .
 \end{aligned}
\end{equation}
Using \eqref{eq:energy}, \(|\mathcal T|\le|A|^3\), and the upper bound for \(|A|\) in Proposition~\ref{th:support}, we obtain
\begin{equation}\label{eq:energy-ratio}
 \begin{aligned}
 \frac{|\mathcal T|^3}
 {\#\{(P_1,P_2,P_3,P_4)\in\mathcal T^4:P_1+P_2=P_3+P_4\}}
 &\le \frac{|A|^9}{N^{3\theta-o(1)}|A|^6}\\
 &\le N^{3(1-\theta)-3\theta+o(1)}
 =N^{6(1/2-\theta)+o(1)} .
 \end{aligned}
\end{equation}

Reiher and Schoen's theorem \cite[Theorem 1.2]{reiher} applies to any nonempty finite subset of an abelian group. Here the group is \(\mathbb Z^2\), the parameter \(K\) is the ratio on the left of \eqref{eq:energy-ratio}, and the error parameter is \(1/4\). The ratio is at least one because any three terms determine the fourth in an additive quadruple. The theorem gives a single nonempty set \(\Omega\subseteq\mathcal T\) satisfying both
\[
 \begin{aligned}
 |\Omega|&\ge\frac34
 \left(\frac{\#\{(P_1,P_2,P_3,P_4)\in\mathcal T^4:P_1+P_2=P_3+P_4\}}
 {|\mathcal T|}\right)^{1/2},\\
 \frac{|\Omega-\Omega|}{|\Omega|}
 &\ll_{1/4}\left(
 \frac{|\mathcal T|^3}{\#\{(P_1,P_2,P_3,P_4)\in\mathcal T^4:P_1+P_2=P_3+P_4\}}
 \right)^4 .
 \end{aligned}
\]
The implicit constant in the second bound depends only on the fixed error parameter \(1/4\).
Substitute \eqref{eq:energy} and \(|\mathcal T|\le|A|^3\) into the bound for \(|\Omega|\), and use \eqref{eq:energy-ratio} in the bound for \(|\Omega-\Omega|/|\Omega|\). We get
\begin{equation}\label{eq:omega-bounds}
 \begin{aligned}
 |\Omega|&\ge \tfrac34
 \left(\frac{N^{3\theta-o(1)}|A|^6}{|A|^3}\right)^{1/2}
 \ge N^{3\theta/2-o(1)}|A|^{3/2},\\
 \frac{|\Omega-\Omega|}{|\Omega|}
 &\ll \bigl(N^{6(1/2-\theta)+o(1)}\bigr)^4
 \le N^{24(1/2-\theta)+o(1)}.
 \end{aligned}
\end{equation}
We will also need a bound for \(|3\Omega-3\Omega|\). For finite nonempty subsets \(U,V\) of an abelian group and positive integers \(j,l\), the Pl\"unnecke--Ruzsa inequality states that \(|jV-lV|\le (|U+V|/|U|)^{j+l}|U|\) \cite[Theorem 1.2]{petridis}. Take \(U=-\Omega,\ V=\Omega\). Then \(|U+V|=|\Omega-\Omega|\) and \(|U|=|\Omega|\). Taking \(j=l=3\) and applying \eqref{eq:omega-bounds}, we obtain
\begin{equation}\label{eq:omega-six}
 \frac{|3\Omega-3\Omega|}{|\Omega|}
 \le\left(\frac{|\Omega-\Omega|}{|\Omega|}\right)^6
 \le N^{144(1/2-\theta)+o(1)} .
\end{equation}

We next select elements of \(A\) that occur in many representations of points in \(\Omega\). Every point of \(\Omega\subseteq\mathcal T\) has at least one such representation, so
\[
 \begin{aligned}
 &\sum_{a\in A}\#\{(u,v)\in A^2:(a+u+v,a^3+u^3+v^3)\in\Omega\}\\
 &=\sum_{P\in\Omega}\#\{(a,u,v)\in A^3:(a+u+v,a^3+u^3+v^3)=P\}
 \ge|\Omega|.
 \end{aligned}
\]
Discard the elements \(a\) with fewer than \(|\Omega|/(2|A|)\) partner pairs. This discards at most \(|A|\cdot|\Omega|/(2|A|)=|\Omega|/2\) triples. At least \(|\Omega|/2\) triples therefore remain.

Choose \((u_0,v_0)\in A^2\) that occurs in the largest number of the remaining triples. Then
\[
 \begin{aligned}
 \frac{|\Omega|}{2}
 &\le \sum_{(u,v)\in A^2}\#\{a\in A:(a,u,v)\text{ remains}\}\\
 &\le |A|^2\max_{(u,v)\in A^2}\#\{a\in A:(a,u,v)\text{ remains}\}\\
 &=|A|^2\,\#\{a\in A:(a,u_0,v_0)\text{ remains}\}.
 \end{aligned}
\]
Let \(S\) be the set of their first coordinates:
\begin{equation}\label{eq:s-definition}
 S=\left\{a\in A:
 \begin{aligned}
 &(a+u_0+v_0,\,a^3+u_0^3+v_0^3)\in\Omega,\\
 &\#\{(u,v)\in A^2:(a+u+v,\,a^3+u^3+v^3)\in\Omega\}
       \ge\frac{|\Omega|}{2|A|}
 \end{aligned}\right\}.
\end{equation}
The maximizing choice of \((u_0,v_0)\), together with \eqref{eq:s-definition}, gives the first inequality below. The remaining bounds use \eqref{eq:omega-bounds} and Proposition~\ref{th:support}.
\begin{equation}\label{eq:slice}
 \begin{aligned}
 |S|&\ge\frac{|\Omega|}{2|A|^2}
 \ge N^{3\theta/2-o(1)}|A|^{-1/2}
 \ge N^{3\theta/2-(1-\theta)/2-o(1)}
 = N^{2\theta-1/2-o(1)},\\
 \frac{|A|}{|S|}&\le
 N^{(1-\theta)-(2\theta-1/2)+o(1)}
 =N^{3(1/2-\theta)+o(1)} .
 \end{aligned}
\end{equation}
By \eqref{eq:s-definition}, every \(a,b\in S\) shares the partner pair \((u_0,v_0)\). Subtracting the corresponding points of \(\Omega\) gives
\begin{equation}\label{eq:displacement}
 \begin{aligned}
 (a-b,a^3-b^3)
 &=(a+u_0+v_0,\,a^3+u_0^3+v_0^3)\\
 &\quad -(b+u_0+v_0,\,b^3+u_0^3+v_0^3)\in\Omega-\Omega .
 \end{aligned}
\end{equation}
\end{proof}

\begin{lemma}\label{le:prime-selection}
Let \(A,\theta\) be as in Lemma~\ref{le:moment-extraction}, and take \(\Omega,S\) satisfying its conclusions. For \(N\) sufficiently large, there exist a prime \(p\in[|S|/16,|S|/8]\) and \(S'\subseteq S\) with \(|S'|\ge3|S|/4\), such that every \(a\in S'\) has at least \(|\Omega|/(4|A|)\) admissible partner pairs.
Call a partner pair \((u,v)\) admissible for \(p\) if
\begin{equation}\label{eq:good-partners}
 p\nmid (u+v),\qquad
 \#\{n\in A:n\equiv w\pmod p\}
 \le \frac{|A|}{|S|}\log^4(2N)
 \quad(w\in\{u,v\}).
\end{equation}

\end{lemma}

\begin{proof}
Let \(p\) run over the primes in \([|S|/16,\,|S|/8]\), and write \(\sum_p1\) for their number. By Lemma~\ref{le:moment-extraction}, \(|S|\ge N^{2\theta-1/2-o(1)}\). Since \(2\theta-1/2>0\), all these primes satisfy \(p\ge|S|/16>3\) once \(N\) is large enough. The prime number theorem states that
\[
 \pi(x):=\#\{q\le x:q\text{ is prime}\}\sim\frac{x}{\log x}
 \qquad(x\to\infty).
\]
Therefore
\begin{equation}\label{eq:prime-count}
 \begin{aligned}
 \sum_p1
 &=\pi(|S|/8)-\pi(|S|/16)+\#\{q:q=|S|/16,\ q\text{ is prime}\}\\
 &=\left(\frac18-\frac1{16}+o(1)\right)\frac{|S|}{\log|S|}
 \sim\frac{|S|}{16\log |S|}
 \gg\frac{|S|}{\log(2N)} \qquad (|S|\le N).
 \end{aligned}
\end{equation}
Let \(j\) be an integer with \(0<|j|\le2N\). The product of its distinct prime divisors divides \(|j|\), and every prime is at least \(2\). Hence
\[
 2^{\#\{p:p\mid j\}}\le\prod_{p\mid j}p\le |j|\le2N.
\]
Taking logarithms and applying \eqref{eq:prime-count}, we obtain
\begin{equation}\label{eq:prime-divisors}
 \#\{p:p\mid j\}\le\log_2(2N)
 \ll\frac{\log^2(2N)}{|S|}\sum_p1 .
\end{equation}
For every \(w\in A\),
\begin{equation}\label{eq:residue-sum}
 \begin{aligned}
 &\sum_p\#\{n\in A:n\equiv w\pmod p\}
 =\sum_{n\in A}\#\{p:p\mid(n-w)\}\\
 &\quad=\underbrace{\sum_p1}_{n=w}
       +\sum_{n\in A\setminus\{w\}}\#\{p:p\mid(n-w)\}\\
 &\quad\overset{\eqref{eq:prime-divisors}}{\ll}
 \left(1+(|A|-1)\frac{\log^2(2N)}{|S|}\right)\sum_p1
 \qquad(0<|n-w|\le N-1\text{ if }n\ne w)\\
 &\quad\ll\frac{|A|}{|S|}\log^2(2N)\sum_p1 \qquad (|S|\le|A|).
 \end{aligned}
\end{equation}
If \(\#\{n\in A:n\equiv w\pmod p\}>(|A|/|S|)\log^4(2N)\), that prime contributes more than \((|A|/|S|)\log^4(2N)\) to the sum in \eqref{eq:residue-sum}. Dividing that sum by the threshold gives
\begin{equation}\label{eq:residue-exceptions}
 \begin{aligned}
 &\#\bigl\{p:\#\{n\in A:n\equiv w\pmod p\}
                  >(|A|/|S|)\log^4(2N)\bigr\}\\
 &\qquad\le
 \frac{\sum_p\#\{n\in A:n\equiv w\pmod p\}}
      {(|A|/|S|)\log^4(2N)}
 \ll \frac{\sum_p1}{\log^2(2N)} .
 \end{aligned}
\end{equation}

A pair \((u,v)\) fails \eqref{eq:good-partners} exactly when at least one of its conditions fails. Thus
\[
 \begin{aligned}
 &\{p:(u,v)\text{ is not admissible for }p\}\\
 &\quad=\{p:p\mid(u+v)\}\\
&\qquad{}\cup\bigcup_{w\in\{u,v\}}
 \bigl\{p:\#\{n\in A:n\equiv w\pmod p\}
               >(|A|/|S|)\log^4(2N)\bigr\}.
 \end{aligned}
\]
Every prime for which \((u,v)\) is not admissible belongs to at least one of the sets on the right. Adding their sizes therefore gives an upper bound for the number of such primes. Since \(2\le u+v\le2N\), \eqref{eq:prime-divisors} gives \(\#\{p:p\mid(u+v)\}\ll(\log^2(2N)/|S|)\sum_p1\). For the residue-class conditions, sum \eqref{eq:residue-exceptions} over \(w\in\{u,v\}\):
\[
 \begin{aligned}
 &\sum_{w\in\{u,v\}}\#\bigl\{p:\#\{n\in A:n\equiv w\pmod p\}
            >(|A|/|S|)\log^4(2N)\bigr\}\\
 &\quad\le \sum_{w\in\{u,v\}}
 \frac{\sum_p\#\{n\in A:n\equiv w\pmod p\}}
      {(|A|/|S|)\log^4(2N)}
 \ll \frac{\#\{u,v\}}{\log^2(2N)}\sum_p1 .
 \end{aligned}
\]
There are at most two values of \(w\), since \(\#\{u,v\}\le2\). Combining the estimates, we get
\begin{equation}\label{eq:pair-exceptions}
 \begin{aligned}
 &\#\{p:(u,v)\text{ is not admissible for }p\}\\
 &\quad\le\#\{p:p\mid (u+v)\}\\
&\qquad{}+\sum_{w\in\{u,v\}}\#\bigl\{p:\#\{n\in A:n\equiv w\pmod p\}
                  >(|A|/|S|)\log^4(2N)\bigr\}\\
 &\quad\ll\left(\frac{\log^2(2N)}{|S|}+\frac{2}{\log^2(2N)}\right)\sum_p1
 =o(1)\sum_p1 .
 \end{aligned}
\end{equation}
The final \(o(1)\) follows from the lower bound for \(|S|\) in Lemma~\ref{le:moment-extraction}. This bound is uniform in \((u,v)\in A^2\).

Fix \(a\in S\). Summing \eqref{eq:pair-exceptions} over its partner pairs and interchanging the two finite sums gives
\begin{equation}\label{eq:partner-exception-sum}
 \begin{aligned}
 &\sum_p\#\{(u,v)\in A^2:(u,v)\text{ partners }a\text{ and is inadmissible for }p\}\\
 &=\sum_{\substack{(u,v)\in A^2\\(a+u+v,a^3+u^3+v^3)\in\Omega}}
         \#\{p:(u,v)\text{ is not admissible for }p\}\\
 &\le o(1)\Bigl(\sum_p1\Bigr)
         \#\{(u,v)\in A^2:(a+u+v,a^3+u^3+v^3)\in\Omega\}.
 \end{aligned}
\end{equation}
If more than half of the partner pairs of \(a\) are inadmissible for a prime \(p\), its contribution exceeds half the total number of partner pairs. Restricting the sum to those primes gives
\[
 \begin{aligned}
 &\tfrac12\,\#\{(u,v)\in A^2:(a+u+v,a^3+u^3+v^3)\in\Omega\}\\
 &\qquad{}\times\#\left\{p:\substack{\text{more than half of the pairs for }a\\\text{are inadmissible for }p}\right\}\\
 &\le \sum_{\substack{p:\,\text{more than half of the pairs for }a\\\text{are inadmissible for }p}}
       \#\left\{(u,v)\in A^2:\substack{(u,v)\text{ partners }a,\\(u,v)\text{ is inadmissible for }p}\right\}\\
 &\le \sum_p\#\{(u,v)\in A^2:(u,v)\text{ partners }a\text{ and is inadmissible for }p\}.
 \end{aligned}
\]
By the hypothesis, \(a\) has at least \(|\Omega|/(2|A|)>0\) partner pairs. Applying \eqref{eq:partner-exception-sum} and dividing by half the total number of partner pairs gives
\begin{equation}\label{eq:half-exceptions}
 \#\{p:\text{more than half of the pairs for }a\text{ are inadmissible for }p\}
 \le o(1)\sum_p1 .
\end{equation}
Summing \eqref{eq:half-exceptions} over \(a\in S\) gives, for sufficiently large \(N\),
\[
 \begin{aligned}
 &\sum_p\#\{a\in S:\text{more than half of its pairs are inadmissible for }p\}\\
 &=\sum_{a\in S}\#\{p:\text{more than half of the pairs for }a\text{ are inadmissible for }p\}
 <\frac{|S|}{4}\sum_p1 .
 \end{aligned}
\]
Consequently, at least one prime \(p\) has fewer than \(|S|/4\) elements \(a\in S\) for which more than half of the partner pairs are inadmissible. Fix this prime. Keep the elements for which at least half of the partner pairs are admissible, and write
\[
 S'=\{a\in S:\text{at least half of its partner pairs are admissible for }p\}.
\]
By the choice of \(p\) and the assumed lower bound for the partner count, we have
\begin{equation}\label{eq:s-prime}
 \begin{aligned}
 &|S'|\ge\frac34|S|,\\
 &\#\{(u,v)\in A^2:(u,v)\text{ is an admissible partner pair of }a\}\ge\frac{|\Omega|}{4|A|}
 \quad(a\in S').
 \end{aligned}
\end{equation}
\end{proof}

\begin{restateproposition}{th:chord-ratio}
Let \(13/27<\theta<1/2\), and let \(N\) be sufficiently large.
If a nonempty set \(A\subseteq\{1,\ldots,N\}\) satisfies \(I_8(\mathbf1_A)\ge N^\theta|A|^4\), then there exist a prime \(p\), a scale \(X=2^j\) with \(j\) a nonnegative integer, and a nonempty set
\[
 E\subseteq\left\{\left(\frac{a-b}{p},\frac{a^3-b^3}{p}\right):
 \begin{array}{l}
 a,b\in A,\ a>b,\ p\mid(a-b),\\
 X\le(a-b)/p<2X
 \end{array}\right\}
\]
such that
\[
 |E|\ge N^{2\theta-1/2-o(1)}.
\]
For every integer \(1\le H\le X\), partition \([X,2X)\) into
\[
 I_j=[X+jH,\,\min\{X+(j+1)H,\,2X\}),\qquad
 j=0,\ldots,\lceil X/H\rceil-1.
\]
If \(h_0\) is the left endpoint of an interval maximizing
\(\#\{(h,k)\in E:h\in I_j\}\), then the nonempty set
\[
 B=\{(h,k)\in E:h_0\le h<h_0+H\}
\]
satisfies
\[
 \frac{|B-B|}{|B|}\le N^{153(1/2-\theta)+o(1)}.
\]
The choices of \(p,X,E\) are independent of \(H\), and the bound is uniform in \(H\) and in the choice of a maximizing interval.
\end{restateproposition}

\begin{proof}[Proof of Proposition \ref{th:chord-ratio}]
Choose \(\Omega,S,(u_0,v_0)\) by Lemma~\ref{le:moment-extraction}, then \(p,S'\) by Lemma~\ref{le:prime-selection}. For this prime \(p\), counting the residue classes of \(S'\) and applying Cauchy--Schwarz gives
\[
 \sum_{j=0}^{p-1}\#\{a\in S':a\equiv j\pmod p\}=|S'|,
 \quad
 \sum_{j=0}^{p-1}\bigl(\#\{a\in S':a\equiv j\pmod p\}\bigr)^2
 \ge\frac{|S'|^2}{p}.
\]
A pair \((a,b)\) satisfies \(p\mid(a-b)\) exactly when both elements lie in the same residue class. Hence
\[
 \#\{(a,b)\in(S')^2:p\mid(a-b)\}
 =\sum_{j=0}^{p-1}\bigl(\#\{a\in S':a\equiv j\pmod p\}\bigr)^2 .
\]
The diagonal pairs are \((a,a)\), one for each \(a\in S'\). Remove these \(|S'|\) pairs. Swapping the coordinates pairs up the remaining cases \(a>b\) and \(a<b\):
\[
 \begin{aligned}
 \{(a,b)\in(S')^2:a>b,\ p\mid(a-b)\}
 &\longrightarrow\{(a,b)\in(S')^2:a<b,\ p\mid(a-b)\},\\
 (a,b)&\longmapsto(b,a).
 \end{aligned}
\]
This map preserves the divisibility condition because \(p\mid(a-b)\iff p\mid(b-a)\), and applying it twice returns the original pair. Thus exactly half of the nondiagonal pairs have \(a>b\). It follows that
\[
 \begin{aligned}
 &\#\{(a,b)\in(S')^2:a>b,\ p\mid(a-b)\}\\
 &=\frac12\left(\sum_{j=0}^{p-1}
    \bigl(\#\{a\in S':a\equiv j\pmod p\}\bigr)^2-|S'|\right)\\
 &\ge\frac12\left(\frac{|S'|^2}{p}-|S'|\right)
 \ge\frac12\left(\frac{(3|S|/4)^2}{|S|/8}-|S|\right)=\frac74|S|.
 \end{aligned}
\]
In the last inequality we used \(p\le|S|/8\) and \(3|S|/4\le|S'|\le|S|\) from \eqref{eq:s-prime} and \(S'\subseteq S\). Lemma~\ref{le:vertical-cubic} says that a chord point determines its endpoints uniquely. Thus the map
\[
 (a,b)\longmapsto
 \left(\frac{a-b}{p},\frac{a^3-b^3}{p}\right)
 \qquad(a,b\in S',\ a>b,\ p\mid(a-b))
\]
is injective.

Since
\[
 \left\{1,\ldots,N\right\}\subseteq
 \bigcup_{j=0}^{\lfloor\log_2N\rfloor}[2^j,2^{j+1}),
\]
the pigeonhole principle gives a scale \(X=2^j\) for which the set
\[
 E=\left\{\left(\frac{a-b}{p},\frac{a^3-b^3}{p}\right):
 \begin{array}{l}
 a,b\in S',\ a>b,\ p\mid(a-b),\\
 X\le(a-b)/p<2X
 \end{array}\right\}
\]
satisfies the first inequality below. The second follows from \eqref{eq:slice}.
\begin{equation}\label{eq:dyadic}
 |E|\ge\frac{7|S|}{4(1+\lfloor\log_2 N\rfloor)}
 \ge N^{2\theta-1/2-o(1)} .
\end{equation}
Keep \(p,X,E\) fixed. We now localize the horizontal coordinate. For any integer \(1\le H\le X\), split \([X,2X)\) into the half-open intervals
\[
 I_j=[X+jH,\,\min\{X+(j+1)H,\,2X\}),\qquad
 j=0,\ldots,\lceil X/H\rceil-1 .
\]
In each interval \(I_j\), we count the points with \(h\in I_j\), namely
\[
 (h,k)=\left(\frac{a-b}{p},\,\frac{a^3-b^3}{p}\right)\in E,
 \qquad a,b\in S',\quad a>b,\quad p\mid(a-b).
\]
Choose an interval \(I_j\) containing the largest number of points of \(E\), and write \(h_0\) for its left endpoint. Since \(h<2X\) for every \((h,k)\in E\), the points in this interval form the set
\[
 B=\{(h,k)\in E:h_0\le h<h_0+H\},\qquad
 |B|=\max_j\#\{(h,k)\in E:h\in I_j\}.
\]
Using
\[
 |E|=\sum_j\#\{(h,k)\in E:h\in I_j\}
 \le \lceil X/H\rceil\,|B|\le(2X/H)|B|
\]
we obtain
\begin{equation}\label{eq:dense-bin}
 h_0\ge X,\qquad |B|\ge \frac{H}{2X}|E|.
\end{equation}

Starting with the empty set, select points \(z=(h,k)\in E\) one at a time. At each step, require the horizontal coordinate to differ by more than \(2H\) from every previously selected coordinate. Stop when no further point can be added. The process stops because \(E\) is finite. Denote the resulting set by \(T\subseteq E\); then
\begin{equation}\label{eq:center-separation}
 |c_1-c'_1|>2H
 \qquad\bigl(c=(c_1,c_2),\ c'=(c'_1,c'_2)\in T,\ c\ne c'\bigr).
\end{equation}
Let \(z=(h,k)\in E\setminus T\). If \(|h-c_1|>2H\) for every \(c=(c_1,c_2)\in T\), we could add \(z\), contrary to the stopping rule. Hence some \(c\in T\) satisfies \(|h-c_1|\le2H\). For \(z\in T\), take \(c=z\). This proves the covering
\[
 E=\bigcup_{c=(c_1,c_2)\in T}
       \{(h,k)\in E:c_1-2H\le h\le c_1+2H\}.
\]

Fix \(c=(c_1,c_2)\in T\). We estimate how many points of \(E\) this center covers. Since \(c_1\in[X,2X)\), the interval \([c_1-2H,c_1+2H]\) meets at least one \(I_j\). Let \(j_-\) and \(j_+\) be the smallest and largest indices of the intervals it meets. They are consecutive:
\[
 I_{j_-},\,I_{j_-+1},\,\ldots,\,I_{j_+}.
\]
Choose \(x_-\in I_{j_-}\cap[c_1-2H,c_1+2H]\) and \(x_+\in I_{j_+}\cap[c_1-2H,c_1+2H]\). For \(j_-<j<j_+\), the interval \(I_j\) is not the possibly shortened final interval, since \(j<\lceil X/H\rceil-1\). Moreover,
\[
 c_1-2H\le x_-<X+jH<X+(j+1)H\le x_+\le c_1+2H .
\]
Hence each interval \(I_j=[X+jH,\,X+(j+1)H)\) with \(j_-<j<j_+\) lies inside \([c_1-2H,c_1+2H]\) and has length \(H\). These intervals are disjoint. If \(j_+-j_-\ge2\), their total length is at most \(4H\):
\[
 (j_+-j_--1)H
 =\sum_{j_-<j<j_+}H
 \le (c_1+2H)-(c_1-2H)=4H .
\]
Including \(I_{j_-}\) and \(I_{j_+}\) gives \(j_+-j_-+1\le6\) intervals in total. If \(j_+-j_-\le1\), there are at most two. Each interval contains at most \(|B|\) points by the choice of \(B\), so
\[
 \begin{aligned}
 \#\{(h,k)\in E:|h-c_1|\le2H\}
 &\le\sum_{j=j_-}^{j_+}\#\{(h,k)\in E:h\in I_j\}\\
 &\le (j_+-j_-+1)|B|\le6|B|.
 \end{aligned}
\]
The sets centered at \(c\in T\) cover \(E\). Summing their bounds and using \eqref{eq:dyadic}, we obtain
\begin{equation}\label{eq:centers}
 |E|\le6|T||B|,\qquad
 |T||B|\ge\frac{|E|}{6}\ge |S|N^{-o(1)} .
\end{equation}

By the definition of \(E\), every \(z=(h,k)\in E\) comes from a pair \(a,b\in S'\). The endpoint uniqueness in Lemma~\ref{le:vertical-cubic} makes this pair unique. Its endpoints satisfy
\[
 a>b,\qquad p\mid(a-b),\qquad
 h=\frac{a-b}{p},\quad k=\frac{a^3-b^3}{p}.
\]
For a center \(c=(c_1,c_2)\in T\subseteq E\), write these endpoints as \(a_c>b_c\). Then
\[
 pc=(pc_1,pc_2)=(a_c-b_c,\,a_c^3-b_c^3).
\]
Since \(S'\subseteq S\), \eqref{eq:s-definition} gives \((x+u_0+v_0,\,x^3+u_0^3+v_0^3)\in\Omega\) for every \(x\in S'\). Use this with \(x=a\) and \(x=b\), and subtract as in \eqref{eq:displacement}:
\[
 \begin{aligned}
 pz=(a-b,\,a^3-b^3)
 &=(a+u_0+v_0,\,a^3+u_0^3+v_0^3)\\
 &\quad -(b+u_0+v_0,\,b^3+u_0^3+v_0^3)\in\Omega-\Omega .
 \end{aligned}
\]
For \(w\in B-B\), choose \(z,z'\in B\subseteq E\) with \(w=z-z'\). Write the endpoints of \(z\) as \(a>b\) and those of \(z'\) as \(a'>b'\). Applying \eqref{eq:displacement} to \(pz\) and \(pz'\), and grouping the positive and negative terms, gives
\begin{equation}\label{eq:chord-lift}
 \begin{aligned}
 pw=pz-pz'
 &=(a+u_0+v_0,\,a^3+u_0^3+v_0^3)
   +(b'+u_0+v_0,\,{b'}^3+u_0^3+v_0^3)\\
 &\quad -(b+u_0+v_0,\,b^3+u_0^3+v_0^3)
   -(a'+u_0+v_0,\,{a'}^3+u_0^3+v_0^3)\\
 &\in 2\Omega-2\Omega .
 \end{aligned}
\end{equation}

Fix a center \(c\in T\) and an admissible partner pair \((u,v)\) of \(a_c\). We group such pairs by their sums and sums of cubes modulo \(p\). Write \(\rho\) for the ordered pair of residue classes
\[
 \begin{gathered}
 \rho_1:=(u+v)+p\mathbb Z,\qquad
 \rho_2:=(u^3+v^3)+p\mathbb Z,\\
 \rho:=(\rho_1,\rho_2)\in(\mathbb Z/p\mathbb Z)^2.
 \end{gathered}
\]
The condition \(\rho_1=0\) is equivalent to \(u+v\equiv0\pmod p\). Admissibility in \eqref{eq:good-partners} gives \(p\nmid(u+v)\), so \(\rho_1\ne0\). Using \((u+v)^3-u^3-v^3=3uv(u+v)\), we get \(3uv\rho_1\equiv\rho_1^3-\rho_2\pmod p\). Since \(p>3\), we can divide by \(3\rho_1\) in \(\mathbb Z/p\mathbb Z\). Thus
\[
 uv\equiv\frac{\rho_1^3-\rho_2}{3\rho_1}\pmod p .
\]
With coefficients in \(\mathbb Z/p\mathbb Z\) and indeterminate \(Z\), this gives
\[
 (Z-u)(Z-v)=Z^2-(u+v)Z+uv
 =Z^2-\rho_1 Z+\frac{\rho_1^3-\rho_2}{3\rho_1}.
\]
Once \(\rho\) is fixed, the quadratic polynomial is fixed. Its roots give at most two ordered pairs of residue classes: two orderings for distinct roots, and one for a repeated root. For each ordered pair of classes, \eqref{eq:good-partners} allows at most \((|A|/|S|)\log^4(2N)\) integer values for each of \(u\) and \(v\). Thus the total number of distinct ordered integer pairs \((u,v)\) with this \(\rho\), across all centers, is at most
\begin{equation}\label{eq:partner-multiplicity}
 2\left(\frac{|A|}{|S|}\log^4(2N)\right)^2 .
\end{equation}

For each \(\rho\) that occurs, choose an integer pair \((u_\rho,v_\rho)\) in these classes that is an admissible partner pair for the largest number of centers. Thus it maximizes
\[
 \#\{c\in T:(u_\rho,v_\rho)\text{ is an admissible partner pair of }a_c\}.
\]
Call a center \(c\) counted in this set compatible with \(\rho\). Before choosing one integer pair for each \(\rho\), \eqref{eq:s-prime} gives the following count of triples \((c,u,v)\):
\begin{equation}\label{eq:partner-total}
 \sum_{c\in T}\#\{(u,v)\in A^2:(u,v)\text{ is an admissible partner pair of }a_c\}
 \ge \frac{|T||\Omega|}{4|A|}.
\end{equation}
Group the triples by \(\rho\), and then by the integer pair \((u,v)\). Within a group with fixed \(\rho,u,v\), the centers \(c\) are distinct. The pair \((u_\rho,v_\rho)\) occurs with at least the average number of centers in its \(\rho\)-group. Divide the triple count in each group by the bound \eqref{eq:partner-multiplicity}, then sum over \(\rho\) and use \eqref{eq:partner-total}. This gives
\begin{equation}\label{eq:compatible}
 \begin{aligned}
 \#\{(\rho,c):c\in T\text{ is compatible with }\rho\}
 &\ge\frac{|T||\Omega|}{4|A|}
 \frac{1}{2\bigl((|A|/|S|)\log^4(2N)\bigr)^2}\\
 &\ge N^{-o(1)}\frac{|T||\Omega||S|^2}{|A|^3}.
 \end{aligned}
\end{equation}

For every compatible \((\rho,c)\) and every \(w\in B-B\), associate the integer point
\begin{equation}\label{eq:injection}
 (u_\rho+v_\rho,\,u_\rho^3+v_\rho^3)+p(c+w).
\end{equation}
Since \(pc=(a_c-b_c,a_c^3-b_c^3)\), we have
\[
 \begin{aligned}
 &(u_\rho+v_\rho,u_\rho^3+v_\rho^3)+pc\\
 &\quad=(a_c+u_\rho+v_\rho,\,a_c^3+u_\rho^3+v_\rho^3)\\
 &\qquad -(b_c+u_0+v_0,\,b_c^3+u_0^3+v_0^3)
            +(u_0+v_0,\,u_0^3+v_0^3).
 \end{aligned}
\]
Compatibility gives the point of \(\Omega\) involving \(a_c\). The inclusion \(b_c\in S'\subseteq S\) gives the point involving \(b_c\), and \eqref{eq:chord-lift} gives the condition on \(pw\):
\[
 \begin{gathered}
 (a_c+u_\rho+v_\rho,a_c^3+u_\rho^3+v_\rho^3)\in\Omega,\\
 (b_c+u_0+v_0,b_c^3+u_0^3+v_0^3)\in\Omega,\qquad
 pw\in2\Omega-2\Omega .
 \end{gathered}
\]
It follows that every point constructed in \eqref{eq:injection} lies in
\[
 (u_0+v_0,\,u_0^3+v_0^3)+(3\Omega-3\Omega).
\]

We claim that the map in \eqref{eq:injection} is injective. Suppose that \((\rho,c,w)\) and \((\rho',c',w')\) give the same point. Reducing modulo \(p\) first determines \(\rho\), after which the chosen integer pair is fixed:
\[
 \begin{aligned}
 &(u_\rho+v_\rho,u_\rho^3+v_\rho^3)+p(c+w)\\
 &\qquad=(u_{\rho'}+v_{\rho'},u_{\rho'}^3+v_{\rho'}^3)+p(c'+w')\\
 &\quad\Longrightarrow\quad
 \rho=\rho'
 \quad\Longrightarrow\quad
 (u_\rho,v_\rho)=(u_{\rho'},v_{\rho'})\\
 &\quad\Longrightarrow\quad c+w=c'+w' .
 \end{aligned}
\]
Since \(w,w'\in B-B\) and \(B\subseteq[h_0,h_0+H]\times\mathbb Z\),
\[
 |c_1-c'_1|=|w'_1-w_1|\le |w'_1|+|w_1|\le 2H .
\]
By \eqref{eq:center-separation}, distinct centers in \(T\) have horizontal coordinates more than \(2H\) apart. Hence \(c=c'\), and then \(w=w'\). This proves injectivity. Since translation does not change the size of a set, we get
\begin{equation}\label{eq:injection-bound}
 \#\{(\rho,c):c\in T\text{ is compatible with }\rho\}\,|B-B|
 \le |3\Omega-3\Omega|.
\end{equation}

Substituting \eqref{eq:compatible} into \eqref{eq:injection-bound} and dividing by \(|B|\) gives the first line below. For the second line, use \eqref{eq:centers} to bound \(|T||B|\). For the last line, use \eqref{eq:omega-six} and \eqref{eq:slice}.
\begin{equation}\label{eq:local-ratio}
 \begin{aligned}
 \frac{|B-B|}{|B|}
 &\le N^{o(1)}\frac{|3\Omega-3\Omega|}{|\Omega|}
             \frac{|A|^3}{|T||B||S|^2}\\
 &\le N^{o(1)}\frac{|3\Omega-3\Omega|}{|\Omega|}
             \left(\frac{|A|}{|S|}\right)^3\\
 &\le N^{(144+3\cdot3)(1/2-\theta)+o(1)}
  =N^{153(1/2-\theta)+o(1)} .
 \end{aligned}
\end{equation}
We chose \(p,X,E\) before choosing \(H\). The argument applies to every integer \(1\le H\le X\) and every interval containing the largest number of points of \(E\), with the same uniform bounds. This proves Proposition~\ref{th:chord-ratio}.
\end{proof}

\begin{restateproposition}{th:chords}
Let \(1/2-1/1550<\theta<1/2\) and let \(N\) be sufficiently large. Suppose that a nonempty set \(A\subseteq\{1,\ldots,N\}\) satisfies \(I_8(\mathbf1_A)\ge N^\theta|A|^4\). Let \(p,X,E\) be as in Proposition~\ref{th:chord-ratio}. Then
\[
 X\ge N^{1/2-155(1/2-\theta)-o(1)}.
\]
For sufficiently large \(N\), the integer
\[
 H=\lfloor XN^{-2/5}\rfloor
\]
satisfies \(1\le H\le XN^{-2/5}\). At this width, the set \(B\) selected in Proposition~\ref{th:chord-ratio} satisfies
\[
 |B|\ge N^{2\theta-9/10-o(1)}.
\]
Its horizontal coordinates lie in \([h_0,h_0+H]\), with \(h_0\ge X\), and
\[
 B\subseteq\left\{\left(\frac{a-b}{p},\frac{a^3-b^3}{p}\right):
 a,b\in A,\ a>b,\ p\mid(a-b)\right\}.
\]
For this same set \(B\), the difference-set estimate in Proposition~\ref{th:chord-ratio} remains valid.
\end{restateproposition}

\begin{proof}[Proof of Proposition \ref{th:chords}]
Since \(1/2-1/1550>13/27\), Proposition~\ref{th:chord-ratio} applies. Fix \(p,X,E\) as in that proposition. We first bound \(X\) from below.
For every integer \(1\le H\le X\), choosing an interval with the largest number of points gives a set \(B\) satisfying \eqref{eq:local-ratio}. Take \(H=X\). There is then just one interval, \([X,2X)\), so \(B=E\) and
\[
 |E-E|/|E|\le N^{153(1/2-\theta)+o(1)} .
\]
There are \(X\) possible horizontal coordinates, and \(\sum_{h=X}^{2X-1}\#\{k:(h,k)\in E\}=|E|\). Cauchy--Schwarz therefore bounds the number of pairs in the same column from below. Group these pairs by their vertical difference \(t\) to obtain
\[
 \begin{aligned}
 \frac{|E|^2}{X}
 &\le \sum_{h=X}^{2X-1}\bigl(\#\{k:(h,k)\in E\}\bigr)^2\\
 &=|E|+
 \sum_{\substack{t\ne0\\(0,t)\in E-E}}
   \#\{(h,k)\in E:(h,k+t)\in E\}\\
 &\le |E|+N^{o(1)}|E-E|.
 \end{aligned}
\]
Here the pairs with equal vertical coordinates contribute \(|E|\). For each nonzero \(t\), Lemma~\ref{le:vertical-cubic}, through \eqref{eq:vertical}, bounds the number of pairs by \(N^{o(1)}\). Explicitly,
\[
 \begin{gathered}
 \#\{((h,k),(h,k'))\in E^2:k'=k\}=|E|,\\
 \#\{(h,k)\in E:(h,k+t)\in E\}\le N^{o(1)}\qquad(t\ne0),\\
 \#\{t\ne0:(0,t)\in E-E\}\le|E-E|.
 \end{gathered}
\]
Divide the column count by \(|E|>0\), then use \eqref{eq:dyadic} and \eqref{eq:local-ratio} with \(B=E\). This gives
\begin{equation}\label{eq:scale-bound}
 \begin{aligned}
 X&\ge \frac{|E|}{1+N^{o(1)}|E-E|/|E|}\\
  &\ge N^{2\theta-1/2-153(1/2-\theta)-o(1)}
   =N^{1/2-155(1/2-\theta)-o(1)} .
 \end{aligned}
\end{equation}
Since \(\theta>1/2-1/1550\),
\[
 \frac12-155\left(\frac12-\theta\right)>\frac12-\frac{155}{1550}=\frac25 .
\]
By \eqref{eq:scale-bound}, \(XN^{-2/5}\ge N^{1/10-155(1/2-\theta)-o(1)}\). The fixed exponent \(1/10-155(1/2-\theta)\) is positive, so \(XN^{-2/5}\to\infty\). We can now take \(H=\lfloor XN^{-2/5}\rfloor\). Once \(N\) is large enough that \(XN^{-2/5}\ge2\),
\[
 1\le \tfrac12 XN^{-2/5}\le XN^{-2/5}-1
 \le H\le XN^{-2/5}\le X .
\]
At this width, choose an interval containing the largest number of points of \(E\), and let \(B\) consist of those points. The bounds \eqref{eq:dense-bin} and \eqref{eq:dyadic} give
\[
 |B|\ge\frac{H}{2X}|E|\ge \tfrac14 N^{-2/5}|E|
       \ge N^{2\theta-1/2-2/5-o(1)}=N^{2\theta-9/10-o(1)} .
\]
For this same \(B\), the estimate \eqref{eq:local-ratio} still holds. Its horizontal coordinates lie in \([h_0,h_0+H]\), with \(h_0\ge X\), and all endpoints lie in \(A\). This proves the proposition.
\end{proof}
\clearpage
\section{Integer chords and the reduction to indicator functions}\label{sec:growth}
\begin{proposition}\label{chord:th:graph}
Fix an integer \(r\ge2\) and a real number \(\kappa\ge1\). Let \(N\ge2\), and let \(B\subseteq\mathbb{Z}^2\) be nonempty, with every point \((h,k)\) satisfying \(1\le h\le N^\kappa,\ |k|\le N^\kappa,\ h\mid k\). Write \(B-B=\{z-z':z,z'\in B\}\). If all horizontal coordinates lie in \([h_0,h_0+H]\), where \(1\le H\le X\le h_0\), then for every \(\varepsilon>0\) there is a constant \(C_{\kappa,r,\varepsilon}\ge1\), depending only on \(\kappa,r,\varepsilon\), such that
\begin{equation}\label{chord:eq:general}
\begin{aligned}
 |B| \le C_{\kappa,r,\varepsilon} N^\varepsilon\bigg(
 &\max_{q\in\mathbb{Z}}\#\{(h,k)\in B:k=qh\}\\
 &+\max_{t\in\mathbb{Z}\setminus\{0\}}\#\{(h,k)\in B:(h,k+t)\in B\}\\
 &+(H/X)^r\max_{(h,k)\in B}\left|\frac{k}{h}\right|\bigg)\left(\frac{|B-B|}{|B|}\right)^{2^r-1}.
\end{aligned}
\end{equation}
\end{proposition}

An ordered additive cube in \(B\) consists of a base point \((h,k)\) and \(r\) directions \((s_i,t_i)\in\mathbb Z^2\) satisfying
\[
  \left(h+\sum_{i\in J}s_i,\ k+\sum_{i\in J}t_i\right)\in B
  \qquad\bigl(J\subseteq\{1,\ldots,r\}\bigr).
\]
We take the empty sum to be zero and allow zero directions and repeated vertices. Once the directions are fixed, call \((h,k)\) admissible if all these vertices lie in \(B\). At each admissible base point, consider
\begin{equation}\label{chord:eq:ell}
 \ell(h,k)=
 \sum_{J\subseteq\{1,\ldots,r\}}(-1)^{r-|J|}
          \frac{k+\sum_{i\in J}t_i}{h+\sum_{i\in J}s_i}.
\end{equation}
We keep the directions fixed whenever we compare \(\ell(h,k)\) and \(\ell(h',k')\).

\begin{lemma}\label{le:integer-levels}
Under the hypotheses of Proposition~\ref{chord:th:graph}, fix the \(r\) directions and use \eqref{chord:eq:ell} on the admissible base points. The values of \(\ell\) are integers, and the number of distinct values is at most
\[
 1+\left\lfloor 2(H/X)^r\max_{(x,y)\in B}|y/x|\right\rfloor.
\]
\end{lemma}

\begin{proof}
Every vertex belongs to \(B\), so its horizontal coordinate divides its vertical coordinate. In particular,
\[
 \frac{k+\sum_{i\in J}t_i}{h+\sum_{i\in J}s_i}\in\mathbb Z
 \quad(\forall J)
 \quad\Longrightarrow\quad \ell(h,k)\in\mathbb Z .
\]

We first bound the variation of \(\ell(h,k)\). Start with the finite geometric series identity
\[
 (1-x)\sum_{j=0}^{r-1}x^j
 =\sum_{j=0}^{r-1}(x^j-x^{j+1})=1-x^r .
\]
Substituting \(x=1-h/(h_0+H/2)\) gives
\[
 1=\frac{h}{h_0+H/2}\sum_{j=0}^{r-1}
       \left(1-\frac h{h_0+H/2}\right)^j
     +\left(1-\frac h{h_0+H/2}\right)^r .
\]
Multiplying both sides by \(k/h\) yields
\begin{equation}\label{chord:eq:slope-decomposition}
 \frac{k}{h}
  =\frac{k}{h_0+H/2}\sum_{j=0}^{r-1}
       \left(1-\frac{h}{h_0+H/2}\right)^j
     +\frac{k}{h}\left(1-\frac{h}{h_0+H/2}\right)^r .
\end{equation}
Since
\[
 h_0\le h\le h_0+H
 \quad\Longrightarrow\quad
 -H/2\le h-(h_0+H/2)\le H/2
\]
we obtain
\[
 \left|1-\frac{h}{h_0+H/2}\right|
  =\frac{|h-(h_0+H/2)|}{h_0+H/2}
  \le\frac{H/2}{h_0+H/2}=\frac{H}{2h_0+H}.
\]
Every vertex has its horizontal coordinate in the same interval, so this bound holds at all vertices.

The first term in the expansion of \(k/h\) is a polynomial in \(h,k\):
\[
 \frac{k}{h_0+H/2}\sum_{j=0}^{r-1}
     \left(1-\frac{h}{h_0+H/2}\right)^j
 =\sum_{j=0}^{r-1}\sum_{a=0}^{j}
       \binom ja\frac{(-1)^a h^a k}{(h_0+H/2)^{a+1}} .
\]
\[
\begin{gathered}
 0\le a\le j\le r-1
 \quad\Longrightarrow\quad \deg(h^ak)=a+1\le r,\\
 a+1=r\ \Longleftrightarrow\ a=j=r-1.
\end{gathered}
\]
For a function \(G\), write \(\Delta_{(s,t)}G(h,k)=G(h+s,k+t)-G(h,k)\). For monomials, the binomial theorem gives
\[
 \Delta_{(s,t)}(h^ak^b)
  =\sum_{\substack{0\le i\le a,\ 0\le j\le b\\(i,j)\ne(a,b)}}
       \binom ai\binom bj h^i k^j s^{a-i}t^{b-j}.
\]
\[
 (i,j)\ne(a,b)\quad\Longrightarrow\quad i+j\le a+b-1
 \quad\Longrightarrow\quad
 \deg\Delta_{(s,t)}G\le\deg G-1
\]
Here the degree comparison is made only for nonzero polynomials. Expanding the differences one at a time gives
\begin{equation}\label{chord:eq:alternating-difference}
 \Delta_{(s_r,t_r)}\cdots\Delta_{(s_1,t_1)}G(h,k)
  =\sum_{J\subseteq\{1,\ldots,r\}}(-1)^{r-|J|}
        G\!\left(h+\sum_{i\in J}s_i,\ k+\sum_{i\in J}t_i\right).
\end{equation}

For \(G\) of class \(C^r\) on a neighborhood of the integration points, the fundamental theorem of calculus expresses each difference as an integral. Applying it \(r\) times and using the chain rule gives
\begin{equation}\label{chord:eq:difference-integral}
\begin{aligned}
 &\Delta_{(s_r,t_r)}\cdots\Delta_{(s_1,t_1)}G(h,k)\\
 &\quad=\int_{[0,1]^r}
 \left[\prod_{i=1}^r(s_i\partial_h+t_i\partial_k)G\right]
 \left(h+\sum_{i=1}^r s_i y_i,\ k+\sum_{i=1}^r t_i y_i\right)
 \,dy_1\cdots dy_r .
\end{aligned}
\end{equation}
Here \(\partial_h,\partial_k\) denote partial differentiation with respect to \(h,k\), respectively. For \(0\le a\le r-1\),
\[
\begin{aligned}
 \prod_{i=1}^r(s_i\partial_h+t_i\partial_k)(h^ak)
 &=
 \left(\prod_{i=1}^r s_i\right)\partial_h^r(h^ak)
 +\sum_{i=1}^r t_i\prod_{\substack{1\le j\le r\\j\ne i}}s_j\,
       \partial_h^{r-1}(h^a)\\
 &=\begin{cases}
 0,&a<r-1,\\
 (r-1)!\displaystyle\sum_{i=1}^r
       t_i\prod_{\substack{1\le j\le r\\j\ne i}}s_j,&a=r-1 .
 \end{cases}
\end{aligned}
\]
By \eqref{chord:eq:difference-integral}, only the term \(h^{r-1}k\) in the polynomial part of \eqref{chord:eq:slope-decomposition} contributes. Using its coefficient from the binomial expansion gives
\[
\begin{aligned}
 &\Delta_{(s_r,t_r)}\cdots\Delta_{(s_1,t_1)}
 \left[\frac{k}{h_0+H/2}\sum_{j=0}^{r-1}
       \left(1-\frac h{h_0+H/2}\right)^j\right]\\
 &\qquad=\frac{(-1)^{r-1}(r-1)!}{(h_0+H/2)^r}
     \sum_{i=1}^r t_i\prod_{\substack{1\le j\le r\\ j\ne i}}s_j .
\end{aligned}
\]
This expression depends only on the directions. Thus in the difference of the two alternating sums \eqref{chord:eq:ell}, the polynomial terms cancel, leaving the remainders from \eqref{chord:eq:slope-decomposition}.

For any vertex \((u,v)\in B\), the remainder satisfies
\begin{equation}\label{chord:eq:remainder}
 \left|\frac vu\left(1-\frac u{h_0+H/2}\right)^r\right|
 \le \left(\frac H{2h_0+H}\right)^r
          \max_{(x,y)\in B}\left|\frac yx\right| .
\end{equation}
The triangle inequality now gives
\[
\begin{aligned}
 |\ell(h',k')-\ell(h,k)|
 &\le \sum_{J\subseteq\{1,\ldots,r\}}
    \left|\frac{k'+\sum_{i\in J}t_i}{h'+\sum_{i\in J}s_i}\right|
    \left|1-\frac{h'+\sum_{i\in J}s_i}{h_0+H/2}\right|^r\\
 &\quad+\sum_{J\subseteq\{1,\ldots,r\}}
    \left|\frac{k+\sum_{i\in J}t_i}{h+\sum_{i\in J}s_i}\right|
    \left|1-\frac{h+\sum_{i\in J}s_i}{h_0+H/2}\right|^r .
\end{aligned}
\]
Apply \eqref{chord:eq:remainder} to each term. There is one term for each subset \(J\), even when some vertices coincide. Since
\[
 \#\{J:J\subseteq\{1,\ldots,r\}\}
    =\sum_{j=0}^r\binom rj=2^r,\qquad
 2h_0+H\ge2X
\]
we obtain
\begin{equation}\label{chord:eq:variation}
\begin{split}
 |\ell(h',k')-\ell(h,k)|
  &\le 2^{r+1}\left(\frac H{2h_0+H}\right)^r\max_{(x,y)\in B}\left|\frac yx\right|\\
  &=2\left(\frac{2H}{2h_0+H}\right)^r\max_{(x,y)\in B}\left|\frac yx\right|\\
  &\le 2(H/X)^r\max_{(x,y)\in B}\left|\frac yx\right| .
\end{split}
\end{equation}

Suppose that at least one admissible base point exists, and list the distinct integer values of \(\ell(h,k)\) as \(\ell_1<\cdots<\ell_m\). Consecutive values differ by at least \(1\), so
\[
 \ell_i\in\mathbb Z,\qquad \ell_{i+1}-\ell_i\ge1,\qquad
 m-1\le\sum_{i=1}^{m-1}(\ell_{i+1}-\ell_i)
   =\ell_m-\ell_1
\]
Combining this with \eqref{chord:eq:variation}, we obtain
\begin{equation}\label{chord:eq:levels}
\begin{aligned}
 &m-1\le 2(H/X)^r\max_{(x,y)\in B}|y/x|,\\
 &m\le1+\left\lfloor2(H/X)^r\max_{(x,y)\in B}|y/x|\right\rfloor,\\
 &2(H/X)^r\max_{(x,y)\in B}|y/x|<1
   \quad\Longrightarrow\quad m=1.
\end{aligned}
\end{equation}
When there is no admissible base point, there are no levels to count. If there is just one level, its value may depend on the directions; it need not be zero.
\end{proof}

\begin{lemma}\label{le:cube-lower}
Let \(B\subseteq\mathbb Z^2\) be finite and nonempty, and write \(C_j(B)\) for the number of ordered \(j\)-dimensional additive cubes in \(B\), including repeated vertices. For each integer \(r\ge1\),
\[
 C_r(B)\ge |B|^{r+1}\left(\frac{|B|}{|B-B|}\right)^{2^r-r-1}.
\]
\end{lemma}

\begin{proof}
Set \(C_0(B)=|B|\).
For fixed \(j\) directions, two admissible base points \(z,z'\) give a \((j+1)\)-dimensional cube by adjoining the direction \(z'-z\). Conversely, every such cube gives these two base points. Thus
\[
\begin{aligned}
 &C_{j+1}(B)\\
 &\quad=\sum_{(s_1,t_1),\ldots,(s_j,t_j)\in B-B}
 \left(\#\left\{z\in B:
 \begin{array}{l}
 z+\sum_{i\in J}(s_i,t_i)\in B\\
 \text{for every }J\subseteq\{1,\ldots,j\}
 \end{array}\right\}\right)^2 .
\end{aligned}
\]
There are at most \(|B-B|^j\) choices for the \(j\) directions. Cauchy--Schwarz gives the first inequality below, and the remaining lines show the induction step:
\[
\begin{aligned}
 &C_{j+1}(B)\ge\frac{C_j(B)^2}{|B-B|^j},\\
 &C_j(B)\ge |B|^{j+1}
          \left(\frac{|B|}{|B-B|}\right)^{2^j-j-1}\\
 &\quad\Longrightarrow\
 C_{j+1}(B)\ge |B|^{j+2}
          \left(\frac{|B|}{|B-B|}\right)^{2(2^j-j-1)+j},\\
 &2(2^j-j-1)+j=2^{j+1}-(j+1)-1 .
\end{aligned}
\]
Starting from \(C_0(B)=|B|\), induction yields
\begin{equation}\label{chord:eq:cube-lower}
 C_r(B)\ge |B|^{r+1}
                \left(\frac{|B|}{|B-B|}\right)^{2^r-r-1}.
\end{equation}
\end{proof}

\begin{proof}[Proof of Proposition~\ref{chord:th:graph}]
Constants \(C\) depend only on the parameters in their subscripts and may increase from line to line. Lemma~\ref{le:integer-levels} gives \eqref{chord:eq:levels} for each tuple of directions. We bound the number of cubes by counting the base points and directions.
We next count cubes for which \(s_1\cdots s_r\ne0\). Reversing the final direction exchanges its two opposite faces:
\[
\begin{aligned}
 &\bigl((h,k);(s_1,t_1),\ldots,(s_r,t_r)\bigr)\\
 &\quad\longleftrightarrow
 \bigl((h+s_r,k+t_r);(s_1,t_1),\ldots,(-s_r,-t_r)\bigr),
 \qquad s_r<0\longleftrightarrow s_r>0 .
\end{aligned}
\]
It is therefore enough to count the cubes with \(s_r>0\) and double the result. Fix the directions and a level \(\ell_0\).

\begin{samepage}
By \eqref{chord:eq:alternating-difference} and \eqref{chord:eq:difference-integral}, taking \(G(h,k)=1/h\) and \(h=x\) gives
\begin{equation}\label{chord:eq:reciprocal}
\begin{split}
 \sum_{J\subseteq\{1,\ldots,r\}}\frac{(-1)^{r-|J|}}{x+\sum_{i\in J}s_i}
 ={}&(-1)^r r!\left(\prod_{i=1}^r s_i\right)\\
 &{}\cdot\int_{[0,1]^r}
      \left(x+\sum_{i=1}^r s_i y_i\right)^{-r-1}
                \,dy_1\cdots dy_r .
\end{split}
\end{equation}
\par\end{samepage}
This identity holds whenever \(x+\sum_i s_i y_i>0\) throughout \([0,1]^r\).
For an admissible base point, the affine function \(h+\sum_i s_i y_i\) attains its minimum on \([0,1]^r\) at a vertex. Hence
\[
 h+\sum_{i=1}^r s_i y_i
 \ge h+\sum_{s_i<0}s_i
 =\min_{J\subseteq\{1,\ldots,r\}}\left(h+\sum_{i\in J}s_i\right)
 \ge h_0>0 .
\]
In \eqref{chord:eq:ell}, the coefficient of \(k\) is the left side of \eqref{chord:eq:reciprocal} with \(x=h\). The integral is positive and \(s_1\cdots s_r\ne0\), so this coefficient is nonzero. We can therefore solve \(\ell(h,k)=\ell_0\) for \(k\):
\[
 \sum_J\frac{(-1)^{r-|J|}}{h+\sum_{i\in J}s_i}\ne0,\qquad
 \ell(h,k)=\ell_0\quad\Longrightarrow\quad k=\frac{P(h)}{Q(h)},
\]
where
\begin{equation}\label{chord:eq:rational}
 \frac{P(x)}{Q(x)}
 =\frac{\displaystyle\ell_0-\sum_{J\subseteq\{1,\ldots,r\}}
           (-1)^{r-|J|}\frac{\sum_{i\in J}t_i}{x+\sum_{i\in J}s_i}}
        {\displaystyle\sum_{J\subseteq\{1,\ldots,r\}}
           \frac{(-1)^{r-|J|}}{x+\sum_{i\in J}s_i}} .
\end{equation}
Reduce this fraction and choose \(P,Q\in\mathbb Z[x]\) coprime in \(\mathbb Q[x]\), with the greatest common divisor of all their coefficients equal to \(1\). If the fraction is identically zero, take \(P=0,Q=1\). At every admissible base point, \(Q(h)\ne0\). The coordinate bounds at the vertices, together with \eqref{chord:eq:ell}, give
\[
 |s_i|,|t_i|\le2N^\kappa,\qquad |\ell_0|\le2^rN^\kappa .
\]
To bound the degrees and coefficients of \(P,Q\), first multiply the numerator and denominator in \eqref{chord:eq:rational} by the product
\[
 \prod_{J\subseteq\{1,\ldots,r\}}\left(x+\sum_{i\in J}s_i\right).
\]
Count repeated factors with their multiplicities. The resulting integer polynomials have degree at most \(2^r\), and their coefficients have absolute value at most \(N^{C_{\kappa,r}}\). To see that reducing the fraction preserves bounds of this form, use Gauss's lemma and the following estimate. If \(G(x)=\sum_{j=0}^m c_jx^j=c_m\prod_{i=1}^m(x-\alpha_i)\) is nonconstant and a nonzero polynomial \(U(x)=\sum u_jx^j\in\mathbb Z[x]\) divides \(G\) in \(\mathbb Z[x]\), then
\[
\begin{aligned}
 &|\alpha_i|\le1+\max_{j<m}|c_j/c_m|,\qquad
   |\text{leading coefficient of }U|\le|c_m|,\\
 &\max_j|u_j|\le
   2^m|c_m|\left(1+\max_i|\alpha_i|\right)^m .
\end{aligned}
\]
It follows that \(P,Q\) have degree at most \(2^r\), with all coefficients bounded in absolute value by \(N^{C_{\kappa,r}}\).

We now distinguish three cases, according to \(Q\) and \(\ell_0\).

Suppose first that \(\deg Q\ge1\). Since \(P,Q\) are coprime, the resultant identity gives \(U,V\in\mathbb Z[x]\) such that
\[
\begin{aligned}
 &UP+VQ=\operatorname{Res}(P,Q)\ne0,\\
 &P(h)=kQ(h)
 \quad\Longrightarrow\quad
 \operatorname{Res}(P,Q)=Q(h)\bigl(kU(h)+V(h)\bigr),\\
 &Q(h)\mid\operatorname{Res}(P,Q),\qquad
 0<|\operatorname{Res}(P,Q)|\le N^{C_{\kappa,r}}.
\end{aligned}
\]
The resultant is the determinant of the Sylvester matrix, so the bounds for the degrees and coefficients give the last inequality. Write \(\tau(m)\) for the number of positive divisors of \(m\ge1\). We use the divisor bound
\[
 \tau(m)\le C_\beta m^\beta\quad(\beta>0).
\]
There are \(2\tau(|\operatorname{Res}(P,Q)|)\) possible nonzero signed divisors \(v\). For each one, \(Q(h)=v\) has at most \(\deg Q\) solutions. Thus
\[
\begin{aligned}
 &\#\{h\in\mathbb Z:Q(h)\ne0,\ Q(h)\mid\operatorname{Res}(P,Q)\}\\
 &\qquad\le2(\deg Q)\tau(|\operatorname{Res}(P,Q)|)
 \le C_{\kappa,r,\varepsilon}N^\varepsilon .
\end{aligned}
\]
Each such \(h\) determines at most one value \(k=P(h)/Q(h)\), giving the same bound for the number of admissible base points.

Suppose next that \(Q\) is a nonzero constant and \(\ell_0=0\). We compare the numerator and denominator of \eqref{chord:eq:rational} as \(x\to+\infty\). By \eqref{chord:eq:reciprocal},
\begin{equation}\label{chord:eq:denominator-limit}
 \lim_{x\to+\infty}x^{r+1}
 \sum_{J\subseteq\{1,\ldots,r\}}
       \frac{(-1)^{r-|J|}}{x+\sum_{i\in J}s_i}
   =(-1)^rr!\prod_{i=1}^r s_i\ne0 .
\end{equation}
For the numerator, first collect the terms containing each \(t_i\):
\[
\begin{aligned}
 &\sum_{J\subseteq\{1,\ldots,r\}}
       (-1)^{r-|J|}\frac{\sum_{i\in J}t_i}{x+\sum_{i\in J}s_i}\\
 &\quad=\sum_{i=1}^r t_i
       \sum_{J\subseteq\{1,\ldots,r\}\setminus\{i\}}
           \frac{(-1)^{r-1-|J|}}{x+s_i+\sum_{j\in J}s_j}.
\end{aligned}
\]
Apply the order-\((r-1)\) form of \eqref{chord:eq:reciprocal} to each inner sum. When \(x>2\sum_i|s_i|\), this gives
\begin{equation}\label{chord:eq:numerator-bound}
 \left|\sum_{J\subseteq\{1,\ldots,r\}}
       (-1)^{r-|J|}\frac{\sum_{i\in J}t_i}{x+\sum_{i\in J}s_i}\right|
 \le \frac{2^r(r-1)!}{x^r}
       \sum_{i=1}^r|t_i|\prod_{\substack{1\le j\le r\\j\ne i}}|s_j| .
\end{equation}
By \eqref{chord:eq:denominator-limit} and \eqref{chord:eq:numerator-bound}, the denominator in \eqref{chord:eq:rational} has order \(x^{-r-1}\) and its numerator is \(O(x^{-r})\) when \(\ell_0=0\). Since \(Q\) is constant, we obtain
\[
 \ell_0=0,\quad Q\text{ is constant}
 \quad\Longrightarrow\quad
 \lim_{x\to+\infty}\frac{P(x)}{x^2}=0
 \quad\Longrightarrow\quad P=0\ \text{or}\ \deg P\le1.
\]
At an admissible base point, \(h\mid k\). Together with \(P(h)=Qk\) from \eqref{chord:eq:rational}, this gives
\[
 h\mid k,\quad P(h)=Qk
 \quad\Longrightarrow\quad h\mid P(h)
 \quad\Longrightarrow\quad h\mid P(0).
\]
If \(P(0)\ne0\), we can count the possible positive values of \(h\) by
\[
 \#\{h\ge1:h\mid P(0)\}
   =\tau(|P(0)|)\le C_{\kappa,r,\varepsilon}N^\varepsilon .
\]
If \(P(0)=0\), the polynomial \(P\) is linear or zero, and
\[
 P(h)=P'(0)h,\qquad k=\frac{P'(0)}Qh,\qquad
 \frac{P'(0)}Q=\frac kh\in\mathbb Z,
\]
All admissible base points then lie on one line of integer slope through the origin. Their number is at most
\[
 \max_{q\in\mathbb{Z}}\#\{(h,k)\in B:k=qh\}.
\]
For fixed directions, this contribution occurs only at the level \(\ell_0=0\).

Finally, suppose that \(Q\) is a nonzero constant and \(\ell_0\ne0\). Using \eqref{chord:eq:denominator-limit} and \eqref{chord:eq:numerator-bound} in \eqref{chord:eq:rational}, with the numerator now tending to \(\ell_0\), gives
\begin{equation}\label{chord:eq:leading-coefficient}
 \deg P=r+1,\qquad
 \frac{P(x)}Q=\sum_{j=0}^{r+1}a_jx^j,\qquad
 a_{r+1}=\frac{(-1)^r\ell_0}{r!\prod_{i=1}^r s_i}\ne0 .
\end{equation}
Use \(P(h)=Qk\) from \eqref{chord:eq:rational} and apply the divisibility condition at each vertex. For every \(J\subseteq\{1,\ldots,r\}\),
\[
\begin{aligned}
 h+\sum_{i\in J}s_i
 &\mid k+\sum_{i\in J}t_i\\
 &\Longrightarrow\
 h+\sum_{i\in J}s_i\mid P(h)+Q\sum_{i\in J}t_i,\\
 P(h)+Q\sum_{i\in J}t_i
 &\equiv P\!\left(-\sum_{i\in J}s_i\right)+Q\sum_{i\in J}t_i
       \pmod{h+\sum_{i\in J}s_i}.
\end{aligned}
\]
Consequently,
\begin{equation}\label{chord:eq:vertex}
 h+\sum_{i\in J}s_i\ \bigm|\ P\!\left(-\sum_{i\in J}s_i\right)+Q\sum_{i\in J}t_i
      \qquad\bigl(J\subseteq\{1,\ldots,r\}\bigr).
\end{equation}
The expression on the right is independent of \(h\): it is the remainder after division by \(x+\sum_{i\in J}s_i\). Its size satisfies
\[
 \left|P\!\left(-\sum_{i\in J}s_i\right)+Q\sum_{i\in J}t_i\right|
 \le N^{C_{\kappa,r}} .
\]
If this remainder is nonzero for some \(J\), \eqref{chord:eq:vertex} makes \(h+\sum_{i\in J}s_i\) a positive divisor of a fixed nonzero integer. Each choice determines \(h\), and hence \(k\). Therefore
\[
 \#\{\text{admissible base points}\}
 \le \tau\!\left(\left|P\!\left(-\sum_{i\in J}s_i\right)
                              +Q\sum_{i\in J}t_i\right|\right)
 \le C_{\kappa,r,\varepsilon}N^\varepsilon .
\]

It remains to consider the case in which all these remainders vanish, with \(Q\) constant and \(\ell_0\ne0\). Then
\begin{equation}\label{chord:eq:zero-remainders}
 \frac{P(-\sum_{i\in J}s_i)}Q=-\sum_{i\in J}t_i
 \qquad\bigl(J\subseteq\{1,\ldots,r\}\bigr).
\end{equation}
In this case we fix the base point and the first \(r-1\) directions, and count the possible final directions.

We shall use a finite-difference identity for polynomials. Take \(m\ge1\) and \(G(x)=b_{m+1}x^{m+1}+b_mx^m+\cdots\). Since
\[
 G^{(m)}(x)=(m+1)!b_{m+1}x+m!b_m,\qquad
 \int_{[0,1]^m}y_i\,dy_1\cdots dy_m=\frac12 ,
\]
applying \eqref{chord:eq:difference-integral} in one variable with steps \(-s_1,\ldots,-s_m\) gives
\begin{equation}\label{chord:eq:poly-face}
\begin{aligned}
 &\sum_{J\subseteq\{1,\ldots,m\}}(-1)^{m-|J|}
       G\!\left(-\sum_{i\in J}s_i\right)\\
 &\quad=(-1)^m\left(\prod_{i=1}^m s_i\right)
       \int_{[0,1]^m}G^{(m)}\!\left(-\sum_{i=1}^m s_i y_i\right)
              \,dy_1\cdots dy_m\\
 &\quad=(-1)^m m!\left(\prod_{i=1}^{m}s_i\right)
     \left(b_m-\frac{m+1}{2}\,b_{m+1}\sum_{i=1}^{m}s_i\right).
\end{aligned}
\end{equation}
To apply \eqref{chord:eq:poly-face} with \(m=r\) and \(G=P/Q\), use \eqref{chord:eq:zero-remainders} and \(r\ge2\) to get
\[
\begin{aligned}
 \sum_J(-1)^{r-|J|}\frac{P(-\sum_{i\in J}s_i)}Q
 &=-\sum_{i=1}^r t_i
       \sum_{J\subseteq\{1,\ldots,r\}\setminus\{i\}}
                   (-1)^{r-1-|J|}\\
 &=-\sum_{i=1}^r t_i(1-1)^{r-1}=0 .
\end{aligned}
\]
Substituting this zero value in \eqref{chord:eq:poly-face} gives
\begin{equation}\label{chord:eq:next-coefficient}
\begin{gathered}
 0=(-1)^rr!\left(\prod_{i=1}^r s_i\right)
     \left(a_r-\frac{r+1}{2}a_{r+1}\sum_{i=1}^r s_i\right)
 ,\\
 \Longrightarrow\quad
 a_r=\frac{r+1}{2}\,a_{r+1}\sum_{i=1}^r s_i .
\end{gathered}
\end{equation}
Apply \eqref{chord:eq:poly-face} a second time with \(m=r-1\) and the polynomial \(G\) below. For an admissible base point and \(J\subseteq\{1,\ldots,r-1\}\), the vertex condition gives \(h+\sum_{i\in J}s_i\ge h_0>0\). Thus the denominator \(-\sum_{i\in J}s_i-h\) is nonzero. Use \eqref{chord:eq:rational} and \eqref{chord:eq:zero-remainders} to evaluate \(G\):
\[
\begin{aligned}
 G(x)&=\frac{P(x)-P(h)}{Q(x-h)}
       =a_{r+1}x^r+(a_r+ha_{r+1})x^{r-1}+\cdots,\\
 G\!\left(-\sum_{i\in J}s_i\right)
   &=\frac{-\sum_{i\in J}t_i-k}{-\sum_{i\in J}s_i-h}
    =\frac{k+\sum_{i\in J}t_i}{h+\sum_{i\in J}s_i}
       \quad\bigl(J\subseteq\{1,\ldots,r-1\}\bigr).
\end{aligned}
\]
Equation \eqref{chord:eq:poly-face} gives
\[
\begin{aligned}
 &\sum_{J\subseteq\{1,\ldots,r-1\}}(-1)^{r-1-|J|}
         \frac{k+\sum_{i\in J}t_i}{h+\sum_{i\in J}s_i}\\
 &\quad=(-1)^{r-1}(r-1)!\left(\prod_{i=1}^{r-1}s_i\right)
 \left(a_r+h a_{r+1}-\frac r2\,a_{r+1}\sum_{i=1}^{r-1}s_i\right).
\end{aligned}
\]
Insert \(a_r\) from \eqref{chord:eq:next-coefficient} and \(a_{r+1}\) from \eqref{chord:eq:leading-coefficient}. After cancellation, we obtain
\begin{equation}\label{chord:eq:face-extension}
\begin{aligned}
 &\ell_0\left(2h+\sum_{i=1}^{r-1}s_i+(r+1)s_r\right)\\
 &\qquad=-2r s_r
 \sum_{J\subseteq\{1,\ldots,r-1\}}(-1)^{r-1-|J|}
         \frac{k+\sum_{i\in J}t_i}{h+\sum_{i\in J}s_i}.
\end{aligned}
\end{equation}

The slopes at all vertices of the base face are integers. The horizontal coordinates of its two opposite vertices give the first bound below, and \eqref{chord:eq:face-extension} gives the second implication:
\[
\begin{aligned}
 2h_0&\le h+\left(h+\sum_{i=1}^{r-1}s_i\right)
           =2h+\sum_{i=1}^{r-1}s_i\le2N^\kappa,\\
 s_r>0,\ \ell_0\ne0
 &\quad\Longrightarrow\quad
 \sum_{J\subseteq\{1,\ldots,r-1\}}(-1)^{r-1-|J|}
         \frac{k+\sum_{i\in J}t_i}{h+\sum_{i\in J}s_i}\ne0 .
\end{aligned}
\]
Multiply \eqref{chord:eq:face-extension} by \(r+1\), and reduce modulo \(2h+\sum_{i=1}^{r-1}s_i+(r+1)s_r\). This yields
\begin{equation}\label{chord:eq:face-divisor}
\begin{aligned}
 &2r\left(2h+\sum_{i=1}^{r-1}s_i\right)
 \sum_{J\subseteq\{1,\ldots,r-1\}}(-1)^{r-1-|J|}
         \frac{k+\sum_{i\in J}t_i}{h+\sum_{i\in J}s_i}\\
 &\qquad\equiv 0 \pmod{2h+\sum_{i=1}^{r-1}s_i+(r+1)s_r}.
\end{aligned}
\end{equation}
The integer on the left depends only on the fixed base point and the first \(r-1\) directions. It is nonzero and satisfies
\[
\begin{aligned}
 0<&\left|2r\left(2h+\sum_{i=1}^{r-1}s_i\right)
 \sum_{J\subseteq\{1,\ldots,r-1\}}(-1)^{r-1-|J|}
         \frac{k+\sum_{i\in J}t_i}{h+\sum_{i\in J}s_i}\right|\\
 &\le 2r\cdot2N^\kappa\cdot2^{r-1}N^\kappa
       =2^{r+1}rN^{2\kappa}.
\end{aligned}
\]
The sum \(2h+\sum_{i=1}^{r-1}s_i\) is at least \(2h_0>0\), and we are counting cubes with \(s_r>0\). Thus \(2h+\sum_{i=1}^{r-1}s_i+(r+1)s_r>0\), and \eqref{chord:eq:face-divisor} makes it a positive divisor of the fixed integer. It increases strictly with \(s_r\), so the divisor bound gives
\[
 \#\{\text{possible values of }s_r\}\le C_{\kappa,r,\varepsilon}N^\varepsilon .
\]
Once \(s_r\) is fixed, \eqref{chord:eq:face-extension} determines \(\ell_0\). To determine \(t_r\), take its coefficient in \eqref{chord:eq:ell} and use the order-\((r-1)\) form of \eqref{chord:eq:reciprocal}:
\[
\begin{aligned}
 \text{coefficient of }t_r
 &=\sum_{J\subseteq\{1,\ldots,r-1\}}
       \frac{(-1)^{r-1-|J|}}{h+s_r+\sum_{i\in J}s_i}\\
 &=(-1)^{r-1}(r-1)!\left(\prod_{i=1}^{r-1}s_i\right)\\
 &\qquad{}\cdot\int_{[0,1]^{r-1}}
         \left(h+s_r+\sum_{i=1}^{r-1}s_i y_i\right)^{-r}
                  \,dy_1\cdots dy_{r-1}\ne0 .
\end{aligned}
\]
This coefficient is nonzero, so \(t_r\) is uniquely determined. There are at most \(|B|^r\) choices for the base point and the first \(r-1\) single-step vertices; these determine the entire base face. Hence
\[
 \#\{\text{ordered base faces}\}\le|B|^r,\qquad
 \#\{\text{cubes in this case}\}
   \le C_{\kappa,r,\varepsilon}N^\varepsilon|B|^r .
\]

We still have to count cubes with a zero horizontal step. If a whole direction is zero, say \((s_i,t_i)=(0,0)\), choose its position, the base point, and the remaining single-step vertices. These choices determine the cube, so
\[
 \#\{\text{cubes with a zero direction}\}\le r|B|^r .
\]
If no direction is zero but some \(s_i=0\), then \(t_i\ne0\). Every admissible base point gives a pair of points with this nonzero vertical difference:
\[
\begin{aligned}
 &t_i\ne0,\qquad (h,k),(h,k+t_i)\in B,\\
 &\#\{\text{admissible base points for these directions}\}\\
 &\qquad\le\#\{(h,k)\in B:(h,k+t_i)\in B\}\\
 &\qquad\le\max_{t\in\mathbb Z\setminus\{0\}}
                 \#\{(h,k)\in B:(h,k+t)\in B\}.
\end{aligned}
\]

For the upper bound, sum the preceding counts over the directions, using \eqref{chord:eq:levels} to bound the number of levels for each direction tuple:
\begin{equation}\label{chord:eq:cube-upper}
\begin{aligned}
 C_r(B)\ \le{}&
 C_{\kappa,r,\varepsilon}N^\varepsilon|B|^r\\
 &+2|B-B|^r\max_{q\in\mathbb{Z}}\#\{(h,k)\in B:k=qh\}\\
 &+|B-B|^r\max_{t\in\mathbb{Z}\setminus\{0\}}\#\{(h,k)\in B:(h,k+t)\in B\}\\
 &+C_{\kappa,r,\varepsilon}N^\varepsilon|B-B|^r
           \left(1+2(H/X)^r\max_{(h,k)\in B}\left|\frac{k}{h}\right|\right).
\end{aligned}
\end{equation}
To combine the terms, use
\[
\begin{gathered}
 \#\{\text{direction tuples}\}\le|B-B|^r,\qquad
 |B|^r\le|B-B|^r,\\
 1\le\max_{q\in\mathbb Z}\#\{(h,k)\in B:k=qh\},\\
 (2^r-r-1)+r=2^r-1 .
\end{gathered}
\]
\begingroup
\emergencystretch=1em
Compare the upper bound \eqref{chord:eq:cube-upper} with the lower bound \eqref{chord:eq:cube-lower}. Divide by \(|B|^r\) and multiply by \((|B-B|/|B|)^{2^r-r-1}\). The resulting inequality is \eqref{chord:eq:general}.\qedhere
\par\endgroup
\end{proof}

\begin{remark}
The three terms in \eqref{chord:eq:general} correspond, respectively, to
\[
\begin{gathered}
 \text{points on a line through the origin}\quad (h,k)\in B,\ k=qh,\\
 \text{a fixed nonzero vertical displacement}\quad (h,k),(h,k+t)\in B,\ t\ne0,\\
 \text{variation of the integer level}\quad
 |\ell(h',k')-\ell(h,k)|
       \le2(H/X)^r\max_{(x,y)\in B}|y/x|.
\end{gathered}
\]
If each vertical column contains at most one point, no nonzero vertical difference occurs, and the second term vanishes. When \(r=2\), \eqref{chord:eq:face-extension} reads
\[
 \ell_0(2h+s_1+3s_2)
   =-4s_2\left(\frac{k+t_1}{h+s_1}-\frac{k}{h}\right).
\]
The expression in parentheses is the difference of the slopes at the two endpoints of the base edge. For \(r\ge3\), the corresponding expression is the alternating difference of order \(r-1\) on the base face.
\end{remark}

\begin{proposition}\label{chord:le:vertical}
Fix integers \(N\ge2,\ p\ge1\), and let \(f\in\mathbb{Z}[x]\) have degree \(d\ge3\). A chord of the curve \((n,f(n))\), directed from \((b,f(b))\) to \((a,f(a))\) and normalized by \(p\), gives the integer point
\[
 (h,k)=\left(\frac{a-b}{p},\frac{f(a)-f(b)}{p}\right),
 \qquad a,b\in\mathbb{Z},\quad 1\le b<a\le N,\quad p\mid a-b .
\]
Let \(B\) be any subset of these normalized chord points. Write \(\tau(m)\) for the number of positive divisors of a positive integer \(m\). For every nonzero integer \(t\),
\[
 \#\{(h,k)\in B:(h,k+t)\in B\}\le 2(d-2)\tau(|t|)^2 .
\]
Consequently, for every \(\eta>0\), this count is at most \(C_{f,\eta}N^\eta\), uniformly in \(t,p,B\).
\end{proposition}

For \(f(x)=x^3\), this is the vertical-difference bound of Lemma~\ref{le:vertical-cubic}.

\begin{proof}
Write \(f(x)=\sum_{j=0}^d c_jx^j\), where \(c_j\in\mathbb Z\) and \(c_d\ne0\). We first count the possible horizontal coordinates. Since \(f\) has integer coefficients,
\[
\begin{aligned}
 &a-b\mid f(a)-f(b)\quad\Longrightarrow\quad h\mid k,\\
 &(h,k),(h,k+t)\in B
 \quad\Longrightarrow\quad h\mid t
 \quad\Longrightarrow\quad
 \#\{\text{possible values of }h\}\le\tau(|t|).
\end{aligned}
\]
For a pair \((h,k),(h,k+t)\), write the smaller endpoints as \(b\) and \(b+w\). Both chords have horizontal increment \(ph\), so their endpoint coordinates are \((b+ph,b)\) and \((b+w+ph,b+w)\). Expanding by the binomial theorem gives
\[
\begin{aligned}
 \frac{f(b+ph)-f(b)}p
 &=\sum_{j=1}^d c_j\sum_{i=1}^j
       \binom ji p^{i-1}h^i b^{j-i}\\
 &=dc_dh\,b^{d-1}+\cdots\in\mathbb Z[b],\\
 \frac{f(b+w+ph)-f(b+w)-f(b+ph)+f(b)}p
 &=t .
\end{aligned}
\]
The expression \((f(b+ph)-f(b))/p\) is an integer polynomial in \(b\). Its values at \(b+w\) and \(b\) differ by a multiple of \(w\). Their difference is \(t\), so
\[
 t\ne0\quad\Longrightarrow\quad
 w\ne0,\quad w\mid t,\quad
 \#\{\text{possible values of }w\}\le2\tau(|t|).
\]
Once \(h,w\) are fixed, the leading term of this difference is
\[
\begin{aligned}
 &\frac{f(b+w+ph)-f(b+w)-f(b+ph)+f(b)}p\\
 &\qquad=d(d-1)c_dhw\,b^{d-2}+\cdots,\qquad
   d(d-1)c_dhw\ne0 .
\end{aligned}
\]
Thus the equation with right-hand side \(t\) has degree exactly \(d-2\) in \(b\), and has at most \(d-2\) solutions. Each choice of \(h,w,b\) determines \(k\), so
\[
\begin{aligned}
 &\#\{(h,k)\in B:(h,k+t)\in B\}\\
 &\le\#\{\text{endpoint data }(h,w,b)\text{ satisfying these conditions}\}\\
 &\le\tau(|t|)\cdot2\tau(|t|)\cdot(d-2)\\
 &=2(d-2)\tau(|t|)^2 .
\end{aligned}
\]
If a pair occurs, its vertical difference satisfies \(1\le|t|\le C_fN^d\). The divisor bound gives
\[
 1\le|t|\le C_fN^d
 \quad\Longrightarrow\quad
 2(d-2)\tau(|t|)^2\le C_{f,\eta}N^\eta .\qedhere
\]
\end{proof}

\begin{remark}
\[
 a-b=ph,\qquad f(a)-f(b)=pk,\qquad
 \frac kh=\frac{f(a)-f(b)}{a-b}.
\]
Points with the same horizontal coordinate \(h\) come from chords with the same horizontal increment \(ph\). We count distinct points \((h,k)\); several of them may lie in one vertical column.

\begin{figure}[htbp]
\centering
\includegraphics[width=.96\linewidth]{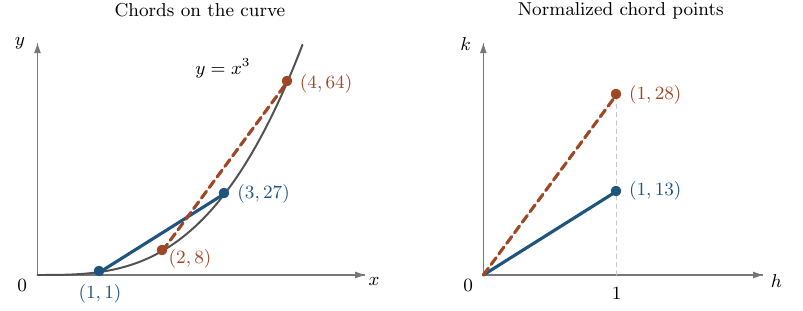}
\caption{For \(d=3\) and \(p=2\), the two chords give the points \((1,13)\) and \((1,28)\). Their difference is \((0,15)\), so they contribute one ordered pair for \(t=15\).}
\label{chord:fig:chord}
\end{figure}
\end{remark}

\begin{proposition}\label{chord:le:slope}
For every integer \(d\ge3\) and every positive integer \(q\), let \(\omega(q)\) denote the number of distinct prime divisors of \(q\). Then
\[
 \#\left\{(a,b)\in\mathbb{Z}^2:1\le b<a,\
                   \frac{a^d-b^d}{a-b}=q\right\}
 \le 2\tau(q)d^{\omega(q)} .
\]
In particular, if \(N\ge2\) is an integer and \(1\le q\le dN^{d-1}\), then for every \(\varepsilon>0\) the left-hand side is at most \(C_{d,\varepsilon}N^\varepsilon\), uniformly in \(N,q\).
\end{proposition}

\begin{proof}
Fix \(q\). We first remove the common factor of \(a,b\):
\[
\begin{gathered}
 g=\gcd(a,b),\qquad a=gu,\quad b=gv,\qquad
 u>v\ge1,\quad\gcd(u,v)=1,\\
 q=g^{d-1}\sum_{j=0}^{d-1}u^{d-1-j}v^j
 \quad\Longrightarrow\quad
 g^{d-1}\mid q,\quad
 \#\{\text{possible values of }g\}\le\tau(q).
\end{gathered}
\]
For each possible \(g\), use the modulus
\[
 m=\frac{q}{g^{d-1}}
  = u^{d-1}+u^{d-2}v+\cdots+v^{d-1}>u^{d-1}.
\]
Since
\[
\begin{aligned}
 &m\equiv u^{d-1}\pmod v,\qquad
 \gcd(v,m)=\gcd(v,u^{d-1})=1,\\
 &u^d-v^d=m(u-v)
 \quad\Longrightarrow\quad(uv^{-1})^d\equiv1\pmod m,
\end{aligned}
\]
where \(v^{-1}\) denotes the multiplicative inverse of \(v\) modulo \(m\). Thus every admissible pair \((u,v)\) gives a root of \(z^d=1\) through the map
\[
 (u,v)\longmapsto uv^{-1}\pmod m
       \in\{z\in(\mathbb Z/m\mathbb Z)^\times:z^d=1\}.
\]
This map is injective. Indeed, if \(uv^{-1}\equiv u'(v')^{-1}\pmod m\), then
\[
\begin{aligned}
 &m\mid uv'-u'v,\qquad
 0<v<u<m^{1/(d-1)},\quad0<v'<u'<m^{1/(d-1)},\\
 &|uv'-u'v|<uu'<m^{2/(d-1)}\le m
       \qquad(d\ge3),\\
 &uv'-u'v=0
 \quad\Longrightarrow\quad \frac uv=\frac{u'}{v'}
 \quad\Longrightarrow\quad (u,v)=(u',v').
\end{aligned}
\]
In the last step, the equal fractions are already reduced, since \(\gcd(u,v)=\gcd(u',v')=1\).

It remains to count these roots. The structure of the unit group modulo a prime power \cite[Chapter 4, Section 1]{ireland-rosen} gives, for a prime \(\pi\) and \(e\ge1\),
\[
\begin{aligned}
 &\#\{z\bmod\pi^e:z^d\equiv1\pmod{\pi^e}\}\\
 &\quad=
 \begin{cases}
 \gcd(d,\pi^{e-1}(\pi-1))\le d,&\pi>2,\\
 1,&\pi=2,\ e=1,\\
 \gcd(d,2)\le2,&\pi=2,\ e=2,\\
 \gcd(d,2)\gcd(d,2^{e-2})\le2d,&\pi=2,\ e\ge3 .
 \end{cases}
\end{aligned}
\]
The Chinese remainder theorem combines the counts for the prime powers:
\[
\begin{aligned}
 \#\{z\bmod m:z^d\equiv1\pmod m\}
 &=\prod_{\pi^e\parallel m}
      \#\{z\bmod\pi^e:z^d\equiv1\pmod{\pi^e}\}\\
 &\le2d^{\omega(m)}\le2d^{\omega(q)} .
\end{aligned}
\]
Here \(\pi^e\parallel m\) means \(\pi^e\mid m\) and \(\pi^{e+1}\nmid m\). There are at most \(\tau(q)\) choices for \(g\), so summing gives
\[
 \#\{(a,b):1\le b<a,\ (a^d-b^d)/(a-b)=q\}
 \le \tau(q)\cdot2d^{\omega(q)} .
\]
Finally, the divisor bound also controls the factor \(d^{\omega(q)}\):
\[
\begin{aligned}
 &2^{\omega(q)}\le\tau(q),\qquad
 d^{\omega(q)}\le\tau(q)^{\log_2d},\\
 &2\tau(q)d^{\omega(q)}
 \le2\tau(q)^{1+\log_2d}
 \le C_{d,\varepsilon}q^{\varepsilon/(d-1)}
 \le C_{d,\varepsilon}N^\varepsilon
       \qquad(q\le dN^{d-1}).
\end{aligned}
\]
\end{proof}

\begin{remark}
For fixed \(p\),
\[
 \#\{(h,k)\in B:k=qh\}
 \le\#\{(a,b):1\le b<a,\ (a^d-b^d)/(a-b)=q\}.
\]
We can now bound both counting terms in Proposition~\ref{chord:th:graph}. Proposition~\ref{chord:le:slope} bounds the number of points on a line through the origin, and Proposition~\ref{chord:le:vertical} bounds the number of pairs with a fixed nonzero vertical difference. Substituting these estimates gives the next result.
\end{remark}

\begin{samepage}
\begin{proposition}\label{chord:th:monomial}
Let \(N\ge2,\ d\ge3,\ r\ge1,\ p\ge1\) be integers, and let \(B\) be any nonempty subset of
\[
 \left\{\left(\frac{a-b}{p},\frac{a^d-b^d}{p}\right):
       a,b\in\mathbb{Z},\ 1\le b<a\le N,\ p\mid a-b\right\}.
\]
If all horizontal coordinates of \(B\) lie in \([h_0,h_0+H]\), with \(1\le H\le X\le h_0\), then for every \(\varepsilon>0\) there is a constant \(C_{d,r,\varepsilon}\ge1\), depending only on \(d,r,\varepsilon\), such that
\begin{equation}\label{chord:eq:local}
 |B|\le C_{d,r,\varepsilon}N^\varepsilon
  \bigl(1+N^{d-1}(H/X)^r\bigr)
  \left(\frac{|B-B|}{|B|}\right)^{2^r-1}.
\end{equation}
In particular, if \(H\le XN^{-(d-1)/r}\), then for every \(\varepsilon>0\),
\begin{equation}\label{chord:eq:expansion}
 |B-B|\ge \frac{1}{C_{d,r,\varepsilon}}N^{-\varepsilon}|B|^{\,1+1/(2^r-1)} .
\end{equation}
\end{proposition}
\end{samepage}

\begin{proof}
For every chord point \((h,k)\in B\),
\begin{equation}\label{chord:eq:slope}
 \frac{k}{h}=\frac{a^d-b^d}{a-b}=\sum_{j=0}^{d-1}a^{d-1-j}b^j .
\end{equation}
\[
 h\mid k,\qquad 1\le h\le N,\qquad |k|\le N^d,\qquad
 1\le k/h\le dN^{d-1}.
\]
Suppose first that \(r\ge2\). The slope range from \eqref{chord:eq:slope} lets us apply Proposition~\ref{chord:le:slope}, while Proposition~\ref{chord:le:vertical} applies to the vertical differences. This gives
\[
\begin{aligned}
 &\max_{q\in\mathbb Z}\#\{(h,k)\in B:k=qh\}\le C_{d,\eta}N^\eta,\\
 &\max_{t\in\mathbb Z\setminus\{0\}}
        \#\{(h,k)\in B:(h,k+t)\in B\}\le C_{d,\eta}N^\eta .
\end{aligned}
\]
In \eqref{chord:eq:general}, take \(\kappa=d\) and power loss \(\varepsilon/2\). Use \(\eta=\varepsilon/2\) in both counts and use \eqref{chord:eq:slope} for the maximum slope. This gives
\[
 \begin{aligned}
 |B|&\le C_{d,r,\varepsilon/2} N^{\varepsilon/2}
       \bigl(2C_{d,\varepsilon/2}N^{\varepsilon/2}
                             +dN^{d-1}(H/X)^r\bigr)
       \left(\frac{|B-B|}{|B|}\right)^{2^r-1}\\
 &\le C_{d,r,\varepsilon}N^\varepsilon\bigl(1+N^{d-1}(H/X)^r\bigr)
       \left(\frac{|B-B|}{|B|}\right)^{2^r-1}.
 \end{aligned}
\]
For \(r=1\), there are at most \(H+1\) integer horizontal coordinates. Apply Cauchy--Schwarz and then use Proposition~\ref{chord:le:vertical} for each nonzero vertical difference:
\[
\begin{aligned}
 \frac{|B|^2}{H+1}
 &\le \sum_{h\in [h_0,h_0+H]\cap\mathbb{Z}}
           \bigl(\#\{k:(h,k)\in B\}\bigr)^2
 \\
 &=|B|+\sum_{\substack{t\ne0\\(0,t)\in B-B}}
          \#\{(h,k)\in B:(h,k+t)\in B\}\\
 &\le |B|+C_{d,\eta}N^\eta|B-B|
 \le C_{d,\eta}N^\eta|B-B|.
\end{aligned}
\]
Also,
\[
 X\le h_0\le N,\qquad
 N^{d-1}H/X\ge H,\qquad H+1\le1+N^{d-1}H/X ,
\]
Divide the preceding count by \(|B|\), multiply by \(H+1\), and take \(\eta=\varepsilon\). The bound for \(H+1\) gives \eqref{chord:eq:local} when \(r=1\).

To obtain the difference-set lower bound, suppose that \(H\le XN^{-(d-1)/r}\). Then
\[
 1+N^{d-1}(H/X)^r\le2 .
\]
Apply \eqref{chord:eq:local} with power loss \((2^r-1)\varepsilon\), and take the \((2^r-1)\)-st root. This gives
\[
\begin{aligned}
 |B|&\le2C_{d,r,(2^r-1)\varepsilon}
     N^{(2^r-1)\varepsilon}
          \left(\frac{|B-B|}{|B|}\right)^{2^r-1},\\
 \frac{|B-B|}{|B|}
 &\ge \left(2C_{d,r,(2^r-1)\varepsilon}\right)^{-1/(2^r-1)}
       N^{-\varepsilon}|B|^{1/(2^r-1)} .
\end{aligned}
\]
Absorb the fixed factor into the constant to obtain \eqref{chord:eq:expansion}.
\end{proof}

\begin{remark}
In \eqref{chord:eq:expansion}, increasing the difference order allows a wider horizontal interval and gives a smaller expansion exponent:
\[
 XN^{-(d-1)/(r+1)}\ge XN^{-(d-1)/r},\qquad
 \frac1{2^{r+1}-1}<\frac1{2^r-1}.
\]
Specializing \eqref{chord:eq:slope} and \eqref{chord:eq:expansion} to \(d=3,\ r=5\) gives
\[
\begin{gathered}
 \max_{(h,k)\in B}k/h\le3N^2,\qquad 2^5-1=31,\\
 H\le XN^{-2/5}
 \quad\Longrightarrow\quad
 |B-B|\ge\frac1{C_{3,5,\varepsilon}}
                N^{-\varepsilon}|B|^{32/31}.
\end{gathered}
\]
For the cubic application, write \(\delta=1/2-\theta\). The derivation of \eqref{eq:scale-bound} uses only Proposition~\ref{th:chord-ratio} and gives
\[
 XN^{-2/r}\ge N^{1/2-2/r-155\delta-o(1)} .
\]
To give this lower bound a positive exponent for \(\delta>0\) sufficiently small, we need \(1/2-2/r>0\), so \(r\ge5\). When the exponent is positive, take \(H=\lfloor XN^{-2/r}\rfloor\). Equations \eqref{eq:dense-bin} and \eqref{eq:dyadic} then give
\[
 |B|\ge N^{1/2-2/r-2\delta-o(1)} .
\]
Compare the lower bound \eqref{chord:eq:expansion} for \(|B-B|\), divided by \(|B|\), with the upper bound \eqref{eq:local-ratio} for \(|B-B|/|B|\). After raising to the \((2^r-1)\)-st power, the comparison yields a contradiction for sufficiently small power losses whenever
\begin{equation}\label{chord:eq:saving-by-order}
 \delta<\frac{1/2-2/r}{153(2^r-1)+2}.
\end{equation}
This condition also ensures \(1/2-2/r-155\delta>0\), since \(153(2^r-1)+2>155\) for \(r\ge5\). For \(r\ge5\), increasing \(r\) by one multiplies the numerator in \eqref{chord:eq:saving-by-order} by
\[
 \frac{1/2-2/(r+1)}{1/2-2/r}
 =1+\frac4{(r+1)(r-4)}\le\frac53 ,
\]
while the denominator more than doubles. Thus the largest threshold in this comparison occurs at \(r=5\). It is \(1/47450\) at \(r=5\) and \(1/57846\) at \(r=6\).
\end{remark}

\begin{remark}
If \(f(x)=\sum_{j=0}^d c_jx^j\in\mathbb Z[x]\), with \(d\ge3\), then
\[
 \left|\frac{f(a)-f(b)}{a-b}\right|
 \le\sum_{j=1}^d|c_j|
             \sum_{i=0}^{j-1}a^{j-1-i}b^i
 \le\left(\sum_{j=1}^d j|c_j|\right)N^{d-1}.
\]
Thus the same argument works for \(f\) if, for every \(\eta>0\) and \(N\ge2\), one has
\[
 \max_{q\in\mathbb{Z}}
 \#\{(a,b)\in\mathbb{Z}^2:1\le b<a\le N,\ f(a)-f(b)=q(a-b)\}
 \le C_{f,\eta}N^\eta .
\]
For \(r\ge2\), combine this assumption with Proposition~\ref{chord:le:vertical} and substitute both estimates in Proposition~\ref{chord:th:graph}. This gives \eqref{chord:eq:local} and \eqref{chord:eq:expansion} for \(f\), with constants that may depend on \(f\).

For \(r=1\), Cauchy--Schwarz over the at most \(H+1\) integer horizontal coordinates gives
\[
\begin{aligned}
 \frac{|B|^2}{H+1}
 &\le \sum_{h\in[h_0,h_0+H]\cap\mathbb Z}
       \bigl(\#\{k:(h,k)\in B\}\bigr)^2\\
 &= |B|+\sum_{\substack{t\ne0\\(0,t)\in B-B}}
       \#\{(h,k)\in B:(h,k+t)\in B\}\\
 &\ll_{f,\varepsilon} N^\varepsilon |B-B|,
\end{aligned}
\]
where the diagonal contributes \(|B|\le |B-B|\) and Proposition~\ref{chord:le:vertical} bounds each summand with \(t\ne0\).
Since \(X\le h_0\le N\), we obtain \eqref{chord:eq:local} in the form
\[
 |B|\ll_{f,\varepsilon} N^\varepsilon(H+1)\frac{|B-B|}{|B|}
 \le N^\varepsilon\left(1+N^{d-1}\frac{H}{X}\right)\frac{|B-B|}{|B|}.
\]

For \(f(n)=\lambda n^d+\mu n+\nu\), where \(\lambda,\mu,\nu\in\mathbb Z\) and \(\lambda\ne0\), we can instead transfer the bounds directly from the monomial case:
\[
\begin{aligned}
 &\left(\frac{a-b}p,\frac{f(a)-f(b)}p\right)
       =(h,\lambda k+\mu h),\\
 &(h,\lambda k+\mu h)=(h',\lambda k'+\mu h')
       \quad\Longrightarrow\quad h=h',\ k=k',\\
 &(h,\lambda k+\mu h)-(h',\lambda k'+\mu h')\\
 &\qquad=(h-h',\lambda(k-k')+\mu(h-h')).
\end{aligned}
\]
The map \((h,k)\mapsto(h,\lambda k+\mu h)\) is additive and injective. It preserves the sizes of both the point set and its difference set, so the same estimates follow.
\end{remark}

\begin{remark}
The restriction \(d\ge3\) is needed. For \(d=2\), Proposition~\ref{chord:th:monomial} fails for every difference order \(r\ge2\). To see this, take \(m\ge1\) and set
\[
\begin{gathered}
 N=16m^2,\quad p=2,\quad h_0=X=N/4,\quad H=m,\\
 a=N/2+h,\quad b=N/2-h,\quad
 N/4\le h\le N/4+m,\quad h\in\mathbb Z,\\
 1\le b<a\le N,\qquad a+b=N,\qquad
 \left(\frac{a-b}{2},\frac{a^2-b^2}{2}\right)=(h,Nh).
\end{gathered}
\]
\[
\begin{gathered}
 B=\{(h,Nh):h\in[X,X+H]\cap\mathbb Z\},\\
 |B|=m+1,\quad |B-B|=2m+1,\\
 H=XN^{-1/2},\qquad N(H/X)^2=1,\qquad
 \frac{|B-B|}{|B|}<2 .
\end{gathered}
\]
If \eqref{chord:eq:local} held for \(d=r=2\), then
\[
 m+1\le 2C_{2,2,\varepsilon}(16m^2)^\varepsilon
             \left(\frac{2m+1}{m+1}\right)^3
       \le 16C_{2,2,\varepsilon}(16m^2)^\varepsilon .
\]
Taking \(\varepsilon=1/4\) gives
\[
 m+1\le32C_{2,2,1/4}\sqrt m
 \quad\Longrightarrow\quad
 \sqrt m\le32C_{2,2,1/4},
\]
which is impossible as \(m\to\infty\). For each fixed \(r>2\), the same example gives
\[
 N(H/X)^r=N^{1-r/2}\le1
 \quad\Longrightarrow\quad
 m+1\le2^{2^r}C_{2,r,\varepsilon}(16m^2)^\varepsilon ,
\]
Taking \(\varepsilon=1/4\) again gives a contradiction as \(m\to\infty\). This example also shows why Proposition~\ref{chord:th:graph} includes the number of points on a line through the origin:
\[
\begin{gathered}
 \max_q\#\{(h,k)\in B:k=qh\}=|B|,\\
 \max_{t\ne0}\#\{(h,k)\in B:(h,k+t)\in B\}=0.
\end{gathered}
\]

Allowing several points in one column gives a counterexample for \(r=1\) as well. Keep \(N,p,X\), but take
\[
\begin{gathered}
 H=1,\qquad a=b+N/2,\qquad b=1,\ldots,N/2,\\
 B=\left\{\left(\frac N4,\frac N4\left(2b+\frac N2\right)\right):
                                     1\le b\le N/2\right\},\\
 |B|=N/2,\qquad |B-B|=2|B|-1,\qquad NH/X=4 .
\end{gathered}
\]
If \eqref{chord:eq:local} held for \(d=2,r=1\), then
\[
 N/2\le10C_{2,1,\varepsilon}N^\varepsilon
 \quad\Longrightarrow\quad
 N^{3/4}\le20C_{2,1,1/4}\quad(\varepsilon=1/4),
\]
which is impossible as \(N\to\infty\). For this set,
\[
 \#\{(h,k)\in B:(h,k+N/2)\in B\}=|B|-1,
\]
Thus the count for a fixed nonzero vertical difference can also be large when \(d=2\). The bound by an arbitrarily small power of \(N\) in Proposition~\ref{chord:le:vertical} requires \(d\ge3\).
\end{remark}

\begin{restateproposition}{th:growth}
Let \(N\ge2\) and \(p\ge1\) be integers, and let \(B\) be any nonempty subset of the following set of integer chord points,
\[
 B\subseteq\left\{\left(\frac{a-b}{p},\,\frac{a^3-b^3}{p}\right):
       a,b\in\mathbb Z,\ 1\le b<a\le N,\ p\mid(a-b)\right\}.
\]
Let \(1\le H\le X\le h_0\). If all horizontal coordinates of \(B\) lie in \([h_0,h_0+H]\) and \(H\le XN^{-2/5}\), then, for every \(\eta>0\), there is a constant \(C_\eta>0\), depending only on \(\eta\), such that
\[
 \frac{|B-B|}{|B|}\ge C_\eta N^{-\eta}|B|^{1/31}.
\]
\end{restateproposition}

\begin{proof}[Proof of Proposition \ref{th:growth}]
We now specialize the monomial estimate to cubic chords. Take \(d=3\) and \(r=5\) in Proposition~\ref{chord:th:monomial}. Then
\[
 2^r-1=31,\qquad
 H\le XN^{-(d-1)/r}=XN^{-2/5}.
\]
Apply \eqref{chord:eq:local} with \(\varepsilon=31\eta\). The condition on \(H\) bounds the factor in parentheses by \(2\), giving
\[
 \begin{aligned}
 |B|
 &\le C_{3,5,31\eta}N^{31\eta}
       \bigl(1+N^2(H/X)^5\bigr)
       \left(\frac{|B-B|}{|B|}\right)^{31}\\
 &\le2C_{3,5,31\eta}N^{31\eta}
       \left(\frac{|B-B|}{|B|}\right)^{31}.
 \end{aligned}
\]
Taking the \(31\)-st root and rearranging, we obtain
\[
 \frac{|B-B|}{|B|}
 \ge(2C_{3,5,31\eta})^{-1/31}N^{-\eta}|B|^{1/31}.
\]
The factor \((2C_{3,5,31\eta})^{-1/31}\) is positive and depends only on \(\eta\). Denoting it by \(C_\eta\) proves the proposition.
\end{proof}

\begin{restateproposition}{prop:weighted}
Fix \(N\ge 1\). If \(D>0\) satisfies \(I_8(\mathbf1_A)\le D|A|^4\) for every \(A\subseteq\{1,\ldots,N\}\), then for all complex coefficients \(b_1,\ldots,b_N\),
\[
 I_8(b)\le D(1+\log N)^4
       \left(\sum_{n=1}^N|b_n|^2\right)^4 .
\]
\end{restateproposition}

\begin{proof}[Proof of Proposition~\ref{prop:weighted}]
First fix a set \(A\) and assign each element a coefficient of modulus \(1\). Expanding the eighth moment and applying \eqref{eq:orthogonality} leaves the same eight-tuples as in \eqref{eq:eight-solutions}. Each coefficient product has absolute value \(1\), so the triangle inequality bounds the moment by \(I_8(\mathbf1_A)\le D|A|^4\).
For general coefficients, express each magnitude \(|b_n|\) as the length of the interval \(0\le t<|b_n|\). Then
\[
 \begin{aligned}
 &\sum_{n=1}^N b_n\exp\!\bigl(2\pi i(nx+n^3y)\bigr)\\
 &\quad=\int_0^\infty
 \sum_{\substack{1\le n\le N\\|b_n|>t}}
 \frac{b_n}{|b_n|}\exp\!\bigl(2\pi i(nx+n^3y)\bigr)\,dt .
 \end{aligned}
\]
The condition \(|b_n|>t\ge0\) excludes \(b_n=0\), so every quotient in this formula is defined.
Apply Minkowski's integral inequality and the bound for coefficients of modulus \(1\):
\[
 I_8(b)^{1/8}\le D^{1/8}
 \int_0^\infty\bigl(\#\{n\in\{1,\ldots,N\}:|b_n|>t\}\bigr)^{1/2}\,dt .
\]
Order the magnitudes by choosing a permutation \(\sigma\) with \(|b_{\sigma(1)}|\ge\cdots\ge|b_{\sigma(N)}|\).
Use \(\sqrt m=\sum_{j=1}^m(\sqrt j-\sqrt{j-1})\) to write the integrand as a finite sum. Interchanging that sum with the integral gives the first equality below. Cauchy--Schwarz then gives
\[
 \begin{aligned}
 &\int_0^\infty\bigl(\#\{n:|b_n|>t\}\bigr)^{1/2}\,dt\\
 &=\sum_{j=1}^N|b_{\sigma(j)}|(\sqrt j-\sqrt{j-1})\\
 &\le\left(\sum_{n=1}^N|b_n|^2\right)^{1/2}
      \left(\sum_{j=1}^N(\sqrt j-\sqrt{j-1})^2\right)^{1/2}\\
 &\le\left(\sum_{n=1}^N|b_n|^2\right)^{1/2}
      \left(\sum_{j=1}^N\frac1j\right)^{1/2}
 \le(1+\log N)^{1/2}\left(\sum_{n=1}^N|b_n|^2\right)^{1/2}.
 \end{aligned}
\]
For the last two inequalities, we used
\[
 (\sqrt j-\sqrt{j-1})^2=(\sqrt j+\sqrt{j-1})^{-2}\le\frac1j,
 \qquad \sum_{j=1}^N\frac1j\le1+\int_1^N\frac{dt}{t}=1+\log N .
\]
Raising to the eighth power proves \eqref{eq:weighted-reduction}.
\end{proof}

\medskip

\subsection*{Acknowledgement} During the preparation of the manuscript,
the authors used GPT for assistance. The authors take responsibility for the mathematical arguments and the final content of the paper.

\bibliographystyle{amsplain}
\bibliography{rf}

\providecommand{\bysame}{\leavevmode\hbox to3em{\hrulefill}\thinspace}
\providecommand{\MR}{\relax\ifhmode\unskip\space\fi MR }
\providecommand{\MRhref}[2]{%
  \href{http://www.ams.org/mathscinet-getitem?mr=#1}{#2}
}
\providecommand{\href}[2]{#2}
\begin{thebibliography}{10}

\bibitem{bourgainkdv}
Jean Bourgain, \emph{Fourier transform restriction phenomena for certain
  lattice subsets and applications to nonlinear evolution equations. {II}. the
  {KdV}-equation}, Geom. Funct. Anal. \textbf{3} (1993), no.~3, 209--262.

\bibitem{bourgain-demeter-guth-vmvt}
Jean Bourgain, Ciprian Demeter, and Larry Guth, \emph{Proof of the main
  conjecture in {Vinogradov's} mean value theorem for degrees higher than
  three}, Ann. of Math. (2) \textbf{184} (2016), no.~2, 633--682.

\bibitem{hardy-wright}
G.~H. Hardy and E.~M. Wright, \emph{An introduction to the theory of numbers},
  sixth ed., Oxford University Press, 2008, Revised by D. R. Heath-Brown and J.
  H. Silverman.

\bibitem{huli}
Yi~Hu and Xiaochun Li, \emph{Discrete {Fourier} restriction associated with
  {KdV} equations}, Anal. PDE \textbf{6} (2013), no.~4, 859--892.

\bibitem{huli-schrodinger}
\bysame, \emph{Discrete {Fourier} restriction associated with {Schr\"odinger}
  equations}, Rev. Mat. Iberoam. \textbf{30} (2014), no.~4, 1281--1300.

\bibitem{huli-fifth}
\bysame, \emph{Local well-posedness of periodic fifth-order {KdV}-type
  equations}, J. Geom. Anal. \textbf{25} (2015), no.~2, 709--739.

\bibitem{hugheswooley}
Kevin Hughes and Trevor~D. Wooley, \emph{Discrete restriction for {$(x,x^3)$}
  and related topics}, Int. Math. Res. Not. IMRN \textbf{2022} (2022), no.~20,
  15612--15631.

\bibitem{ireland-rosen}
Kenneth Ireland and Michael Rosen, \emph{A classical introduction to modern
  number theory}, second ed., Graduate Texts in Mathematics, vol.~84, Springer,
  New York, 1990.

\bibitem{laiding}
Xudong Lai and Yong Ding, \emph{A note on the discrete {Fourier} restriction
  problem}, Proc. Amer. Math. Soc. \textbf{146} (2018), no.~9, 3839--3846.

\bibitem{li-mean-value}
Xiaochun Li, \emph{A remark on the mean value inequalities related to
  {Waring's} problem}, Preprint, 2023.

\bibitem{li-stein-decoupling}
\bysame, \emph{A {Stein--Tomas} type estimate and a decoupling inequality},
  Preprint, 2023.

\bibitem{li-stein-2025}
\bysame, \emph{A remark on {Stein--Tomas} type restriction theorems}, Acta
  Math. Sin. (Engl. Ser.) \textbf{41} (2025), no.~1, 122--130.

\bibitem{petridis}
Giorgis Petridis, \emph{New proofs of {Pl{\"u}nnecke}-type estimates for
  product sets in groups}, Combinatorica \textbf{32} (2012), no.~6, 721--733.

\bibitem{reiher}
Christian Reiher and Tomasz Schoen, \emph{Note on the theorem of {Balog},
  {Szemer\'edi}, and {Gowers}}, Combinatorica \textbf{44} (2024), no.~3,
  691--698.

\bibitem{wooley}
Trevor~D. Wooley, \emph{Mean value estimates for odd cubic {Weyl} sums}, Bull.
  Lond. Math. Soc. \textbf{47} (2015), no.~6, 946--957.

\end{thebibliography}
\end{document}